\documentclass{amsart}
\usepackage{enumerate,enumitem}
\usepackage{amsthm, amssymb, amsmath}
\usepackage[initials]{amsrefs}
\usepackage{graphicx, tikz}
\usetikzlibrary{patterns}
\usepackage{mathrsfs}
\usepackage{mathtools}
\usepackage{subfig}
\usepackage{tabularx}
\usepackage{fancyhdr}
\usepackage{xcolor}
\usepackage{hyperref}

\hypersetup{
    colorlinks,
    linkcolor={red!75!black},
    citecolor={blue!70!black},
    urlcolor={blue!80!black}
}

\newtheorem{theorem}{Theorem}[section]
\newtheorem{prop}[theorem]{Proposition}
\newtheorem{cor}[theorem]{Corollary}
\newtheorem{lm}[theorem]{Lemma}
\newtheorem{manualtheoreminner}{Theorem}
\newenvironment{manualtheorem}[1]{%
  \IfBlankTF{#1}
    {\renewcommand{\themanualtheoreminner}{\unskip}}
    {\renewcommand\themanualtheoreminner{#1}}%
  \manualtheoreminner
}{\endmanualtheoreminner}
\theoremstyle{definition}
\newtheorem{defin}[theorem]{Definition}
\newtheorem*{rem*}{Remark}
\newtheorem{rem}[theorem]{Remark}

\newcommand{\D}{\mathbb{D}}
\newcommand{\C}{\mathbb{C}}
\newcommand{\R}{\mathbb{R}}
\newcommand{\N}{\mathbb{N}}
\renewcommand{\H}{\mathbb{H}}

\renewcommand{\epsilon}{\varepsilon}
\renewcommand{\subset}{\subseteq}

\numberwithin{equation}{section}

\title[Semigroup models for discrete iteration]{Semigroup models for discrete holomorphic iteration in the unit disc}
\author[A. Christodoulou]{Argyrios Christodoulou$^1$}
\address{Department of Mathematics, Aristotle University of Thessaloniki, 54124, Thessaloniki, Greece}
\email{argyriac@math.auth.gr}
\author[K. Zarvalis]{Konstantinos Zarvalis}
\address{Department of Mathematics, Aristotle University of Thessaloniki, 54124, Thessaloniki, Greece}
\email{zarkonath@math.auth.gr}
\thanks{$^1$ Corresponding author}

\begin{document}
	
	\begin{abstract}
		The main goal of this article is to develop a technique that allows us to partially embed the orbit of any holomorphic self-map $f$ of the disc, into a semigroup which captures the asymptotic behaviour of the orbit. This extends the semigroup-fication procedure introduced by Bracci and Roth to non-univalent functions. We use our technique in order to obtain sharp estimates for the rate with which the orbits of $f$ converge to the attracting fixed point; a fundamental, yet underdeveloped, concept in discrete iteration. Moreover, our results allow us to evaluate the slope of the orbits of $f$, and prove that they behave similarly to quasi-geodesic curves precisely when they converge non-tangentially. 
	\end{abstract}
	
	\maketitle

    \tableofcontents
    \addtocontents{toc}{\protect\setcounter{tocdepth}{1}}

\section{Introduction}
One of the most prominent results on the topic of holomorphic iteration in the unit disc $\D$ of the complex plane is the Denjoy--Wolff Theorem, which states that the iterates of a holomorphic self-map $f$ of $\D$ converge to a unique point $\tau\in\overline{\D}$, whenever $f$ is not conjugate to a Euclidean rotation. This result, however, does not provide any information on the manner in which the iterates approach the \emph{Denjoy--Wolff point} $\tau$. Deciphering the precise nature of this convergence has been the topic of research for several decades \cites{Arosio-Bracci, Baker-Pommerenke, BMS, Bracci-Poggi, CCZRP, parabolic-zoo, CDP2, Fran, Pommerenke-Iteration}; yet many of its elements remain unclear.

On the other hand, in the theory of semigroups of holomorphic functions---another branch of holomorphic dynamics---the asymptotic behaviour of the trajectories of a semigroup is very well-understood. This is the culmination of almost five decades of research and numerous articles, such as \cites{BP,Betsakos-Asymptotic, BCDM-Rates,Bracci-Speeds, BCDG, BCDMGZ, CDM, Kelgiannis}, to name a few.

Studying the connections between discrete and continuous dynamics has a long history, dating back to Koenigs, and has been an ongoing process for several decades (see, for example, \cite{BCDM, Bracci-Roth, Cowen, KRS, SZ, CDMG} and references therein). This article contributes to this endeavour by developing a technique that extends a result of Bracci and Roth \cite{Bracci-Roth}, which allows us to partially embed the orbits of any holomorphic self-map of $\D$ into a trajectory of a semigroup of holomorphic functions.

Our construction enables us to draw from the large pool of results and techniques present in the theory of semigroups in order to evaluate the slope and the rate at which the iterates of the self-map approach the Denjoy--Wolff point; two fundamental concepts in iteration theory. Outside the realm of iteration theory, we use this technique to study the hyperbolic geometry of a certain class of multiply connected domains, and to obtain sharp estimates on the growth rate of norms of composition operators in classical function spaces.

To formally state our results, we start by defining the \emph{iterates} of a holomorphic function $f\colon \D\to\D$ as the $n$-fold compositions $f^n\vcentcolon=f\circ f \circ\cdots\circ f$, for $n\in\N$. We also write $f^0=\mathrm{Id}_\D$. By the Denjoy--Wolff Theorem, if $f$ does not have any fixed points in $\D$, there exists a unique $\tau\in\partial\D$ such that $\{f^n(z)\}$ converges to $\tau$ for all $z\in\D$. We say that such a holomorphic map $f$ is \emph{non-elliptic} and the point $\tau$ is called its \emph{Denjoy--Wolff point}.

An important tool in iteration theory of the unit disc is the ``linearisation" of non-elliptic maps, described as follows. A domain $\Omega\subset\C$ is called \emph{starlike at infinity} if $\Omega+t\subset\Omega$ for all $t\geq0$. Similarly, $\Omega$ is called \emph{asymptotically starlike at infinity} if $\Omega+1\subset\Omega$ and the domain $\widetilde{\Omega}\vcentcolon=\bigcup_{n=1}^\infty\left(\Omega-n\right)$ is starlike at infinity. For any non-elliptic $f\colon\D\to\D$, there exist a domain $\Omega$ asymptotically starlike at infinity and an onto holomorphic map $h\colon\D\to\Omega$ so that $h\circ f =h +1$, called a \textit{Koenigs domain} and a \textit{Koenigs function} for $f$, respectively. Both $\Omega$ and $h$ are unique up to translation, and there are essentially only three possibilities for the domain $\widetilde{\Omega}$, that determine the \emph{type} of $f$: if $\widetilde{\Omega}$ is a horizontal strip, $f$ is called \emph{hyperbolic}; if $\widetilde{\Omega}$ is a horizontal half-plane, $f$ is called \emph{parabolic of positive hyperbolic step}; and if $\widetilde{\Omega}$ is the complex plane, $f$ is called \emph{parabolic of zero hyperbolic step}.

The theory surrounding the Koenigs domain and Koenigs function is the product of the work of Valiron \cite{Valiron}, Pommerenke \cite{Pommerenke-Iteration}, Baker and Pommerenke \cite{Baker-Pommerenke} and Cowen \cite{Cowen} (see also \cite{Arosio-Bracci}). The article \cite{Cowen}, in particular, proves the existence of domains on which a self-map of the disc is well-behaved, that are key to our analysis. A simply connected domain $U\subset \D$ is called \emph{a fundamental domain} of a holomorphic map $f\colon \D\to\D$ if $f$ is univalent on $U$, $f(U)\subset U$, and $\bigcup_{n\in\N}f^{-n}(U)=\D$, where $f^{-n}(U)$ denotes the preimage of $U$ under the iterate $f^n$.

Our first step is to find a fundamental domain of $f$ that interacts particularly well with a Koenigs function of a non-elliptic $f$, and whose hyperbolic geometry is comparable with the hyperbolic geometry of $\D$ close to the Denjoy--Wolff point. To describe this, we equip a domain $D$, whose complement $\C\setminus D$ contains at least two points, with the hyperbolic distance $d_D(\cdot,\cdot)$ induced by the hyperbolic metric $\lambda_D(z)\lvert dz\rvert$ (see Section \ref{sect: hyperbolic geometry}). 

\begin{manualtheorem}{A}\label{thm:main A}
    Let $f\colon \D\to\D$ be a non-elliptic map with Denjoy--Wolff point $\tau\in\partial\D$, Koenigs function $h$ and Koenigs domain $\Omega$. There exists a fundamental domain $V\subset \D$ of $f$ such that $h$ is univalent on $V$, $h(V)$ is starlike at infinity and 
    \begin{equation}\label{eq: semigroupfication thm 1 eq1}
    \bigcup_{t\geq0}\left( h(V)-t\right)=\bigcup_{n\in\N}\left( \Omega-n\right).
    \end{equation}
    When $f$ is hyperbolic or parabolic of zero hyperbolic step, any sequence $\{z_n\}\subset\D$ converging non-tangentially to $\tau$ is eventually contained in $V$, and
    \begin{equation}\label{eq: semigroupfication thm1 eq2}
        \lim_{n\to+\infty}\frac{\lambda_\D(z_n)}{\lambda_V(z_n)}=1.
    \end{equation}
\end{manualtheorem}

Using Theorem \ref{thm:main A} we can define 
\begin{equation}\label{eq:semigroup-fication def}
    \phi_t(z)\vcentcolon=h\lvert_V^{-1}\left(h\lvert_V(z)+t\right), \quad \text{for any } z\in V \text{ and all } t\geq0.
\end{equation}
Since $h$ is univalent on $V$ and $h(V)$ is starlike at infinity, $(\phi_t)$ is a well-defined semigroup of holomorphic functions in $V$ (see Section \ref{sect:semigroups}). Thus, the restriction of $f$ in $V$ can be embedded into the semigroup $(\phi_t)$, which we call the \emph{semigroup-fication of $f$ in $V$}. 

This terminology is borrowed from Bracci and Roth \cite{Bracci-Roth} who constructed the semigroup-fication of univalent self-maps of $\D$. In fact, for the special case of a univalent, non-elliptic $f$, our procedure recovers the semigroup-fication of \cite{Bracci-Roth} (see Section \ref{section:fundamental domains and semigroupfication}).

Since the linearisation and the type of a semigroup are defined similarly to the case of self-maps, we can use \eqref{eq: semigroupfication thm 1 eq1} to show that $f$ and its semigroup-fication $(\phi_t)$ have the same type. Moreover, the limit \eqref{eq: semigroupfication thm1 eq2} in Theorem \ref{thm:main A} allows us to show that, in many cases, the hyperbolic distances $d_\D$ and $d_V$ are Lipschitz equivalent close to the Denjoy--Wolff point of $f$. This implies that the sequence $\{f^n(z)\}$ and the curve $\phi_t(z)$, with $t\in[0,+\infty)$, exhibit similar asymptotic behaviour in $\D$, for all $z\in V$. So, even though the semigroup $(\phi_t)$ is only defined on a subdomain of $\D$, it ``captures" the dynamical properties of $f$.  

These core elements of the semigroup-fication of $f$ in $V$ are collected in the following theorem.

\begin{manualtheorem}{B}\label{thm:main B}
    Let $f\colon\D\to\D$ be a non-elliptic map with Denjoy--Wolff point $\tau\in\partial\D$, and let $(\phi_t)$ be the semigroup-fication of $f$ in $V$, where $V$ is the fundamental domain provided by Theorem \ref{thm:main A} and the semigroup $(\phi_t)$ is the one given in \eqref{eq:semigroup-fication def}. Then
    \begin{enumerate}[label=\textup{(\alph*)}]
        \item $\phi_1\equiv f\lvert_V$ and $\phi_n(z)=f^n(z)$, for any $z\in V$ and all $n\in\N$;
        \item $f$ and $(\phi_t)$ have the same type (hyperbolic, parabolic of zero hyperbolic step or parabolic of positive hyperbolic step);
        \item $\lim\limits_{t\to+\infty} \lvert \phi_t(z)-\tau\rvert=0$, for all $z\in V$; and
        \item for any $z\in V$, $\{f^n(z)\}$ converges to $\tau$ non-tangentially if and only if $\phi_t(z)$ converges to $\tau$ non-tangentially as $t\to+\infty$.
    \end{enumerate}
\end{manualtheorem}

We now describe how the extensive literature in the theory of semigroups can be used to examine the manner in which orbits approach the Denjoy--Wolff point.

First, given $\{z_n\}\subset\D$ converging to $\zeta\in\partial\D$, we define the slope of $\{z_n\}$ as the set of accumulation points of the sequence $\left\{\mathrm{arg}(1-\overline{\zeta}z_n)\right\}$, which we denote by $\mathrm{Slope}_\D(z_n)\subset [-\tfrac{\pi}{2},\tfrac{\pi}{2}]$. We say that a curve $\gamma\colon[0,+\infty)\to\D$ \emph{lands} at a point $\zeta\in\partial\D$ if $\lim_{t\to+\infty}\gamma(t)=\zeta$. The \emph{slope} of $\gamma$ is defined as the cluster set of $\left\{\mathrm{arg}(1-\overline{\zeta}\gamma(t))\colon t\geq0\right\}$, as $t\to +\infty$, and is denoted by $\mathrm{Slope}_\D(\gamma)$. 

Also, a curve $\gamma\colon[0,+\infty)\to\D$ is called a \emph{hyperbolic quasi-geodesic} of $\D$ if there exist $A\geq 1$ and $B\geq 0$ so that
\begin{equation}\label{eq:quasi-geodesic}
\ell_\D(\gamma;[t_1,t_2])\leq Ad_\D(\gamma(t_1),\gamma(t_2))+B,\quad\text{for all}\ 0\leq t_1<t_2,
\end{equation}
where $\ell_\D(\gamma;[t_1,t_2])$ denotes the hyperbolic length of $\gamma$ between $\gamma(t_1)$ and $\gamma(t_2)$. The concept of quasi-geodesic curves originates in Gromov's hyperbolicity theory (see, for example, \cite{Gromov}), and they constitute a class of curves closely related to---yet far more wieldy than---the geodesics of a metric space. Recently, quasi-geodesics were employed in holomorphic dynamics \cites{BCDMGZ,Z}, in order to study the asymptotic behaviour of semigroups of holomorphic functions in $\D$.

Our first application of Theorem \ref{thm:main B} shows that the slope of any orbit of a non-elliptic $f$ is completely determined by the slope of the trajectories of its semigroup-fication. In particular, this demonstrates that the orbits of $f$ can be embedded into $f$-invariant, Lipschitz curves of the same slope. 

\begin{manualtheorem}{C}\label{thm:main slope}
    Let $f\colon\D\to\D$ be a non-elliptic map with Denjoy--Wolff point $\tau\in\partial\D$, and $(\phi_t)$ its semigroup-fication in $V$, where $V$ is the fundamental domain provided by Theorem \ref{thm:main A} and the semigroup $(\phi_t)$ is the one given in \eqref{eq:semigroup-fication def}. For any $z\in\D$, there exists some $n_0\in\N$ such that $\eta_z\colon [0,+\infty)\to \D$ with $\eta_z(t)=\phi_t(f^{n_0}(z))$ is a well-defined, Lipschitz curve that lands at $\tau$ and satisfies:
    \begin{enumerate}[label=\rm (\alph*)]
        \item $f^n(z)=\eta_z(n-n_0)$, for all $n\geq n_0$;
        \item $f(\eta_z([0,+\infty)))\subset \eta_z([0,+\infty))$; and
        \item $\mathrm{Slope}_\D(f^n(z))=\mathrm{Slope}_\D(\eta_z)$.
    \end{enumerate}
    Moreover, $\eta_z$ is a hyperbolic quasi-geodesic of $\D$ if and only if $\{f^n(z)\}$ converges to $\tau$ non-tangentially.
\end{manualtheorem}

Theorem \ref{thm:main slope} also tells us that the orbits of $f$ can be embedded into $f$-invariant quasi-geodesics of $\D$, whenever they converge non-tangentially. Using this we prove that, in this case, the sum of the hyperbolic distances between consecutive terms of the orbit is controlled by the distance between the starting and the ending term. This property can be thought of as a discrete analogue of the fact that non-tangential trajectories of semigroups are quasi-geodesic curves (see \cite[Theorem 1.2]{BCDMGZ}). 

\begin{cor}\label{coro:main slope coro}
For any non-elliptic map $f:\D\to\D$, the following conditions are equivalent:
    \begin{enumerate}[label=\rm (\alph*)]
    \item For any $z\in\D$, there exist constants $A\geq1$ and $B\geq0$ so that for all integers $0\leq n<m$, we have 
    \begin{equation} \label{eq:quasi-lemma eq}
    \sum_{k=n}^{m-1}d_{\D}(f^k(z),f^{k+1}(z)) \leq Ad_{\D}(f^n(z),f^m(z))+B.
    \end{equation}
    \item The orbit $\{f^n(z)\}$ converges to the Denjoy--Wolff point of $f$ non-tangentially, for some $z\in\D$.
    \end{enumerate}
\end{cor}

Next, we turn our attention to the rate with which orbits move towards the Denjoy--Wolff point. Using established results on the rates of convergence of semigroups we obtain the following estimate.

\begin{manualtheorem}{D}\label{thm:main rates}
    Let $f:\D\to\D$ be a non-elliptic map whose Koenigs domain is not the whole complex plane. For every $z\in\D$ and every $\epsilon>0$, there exists a constant $c\vcentcolon=c(z,\epsilon)$ such that
    \begin{equation}\label{eq:hyperbolic rate eq}
    d_{\D}(z,f^n(z))\ge \dfrac{1}{4+\epsilon}\log n + c, \quad\text{for all }n\in\N.
    \end{equation}
\end{manualtheorem}

The inequality in Theorem \ref{thm:main rates} is best possible, in the sense that there exists a non-elliptic map $f$ for which 
\[
\lim_{n\to+\infty}\frac{d_{\D}(z,f^n(z))}{\log n}=\frac{1}{4}, 
\]
(see Remark \ref{rem:sharpness of rates}). Moreover, it is currently not known whether there exist non-elliptic maps whose Koenigs domain is $\C$. If they do exist, then one would probably need to employ techniques different from ours in order to obtain a result similar to Theorem \ref{thm:main rates}.

It is known that when $f$ is hyperbolic or parabolic of positive hyperbolic step, the estimate in Theorem \ref{thm:main rates} can be improved; i.e. the term  $\frac{1}{4+\epsilon}\log n$ may be replaced by a larger quantity (see Section \ref{section:convergence rates}). The behaviour of parabolic maps of zero hyperbolic step, however, is notoriously chaotic and no general estimate on their rate of convergence exists in the literature; especially for non-univalent maps. This is the main contribution of Theorem \ref{thm:main rates}.

The quantity $d_\D(z,f^n(z))$ is sometimes called the \emph{divergence rate} of $f$ since it measures how quickly $f^n(z)$ moves away from $z$ (see \cite{Arosio-Bracci} or \cite[Section 9.1]{BCDM-Book}). The term \emph{rate of convergence} is typically reserved for Euclidean quantities such as the following, which is an equivalent form of Theorem \ref{thm:main rates}.

\begin{manualtheorem}{D*}\label{thm:main rates 2}
    Let $f:\D\to\D$ be a non-elliptic map whose Koenigs domain is not the whole complex plane. For every $z\in\D$ and every $\epsilon>0$, there exists a positive constant $c\vcentcolon=c(z,\epsilon)$ such that
    \[
    1-\lvert f^n(z)\rvert \leq c\ n^{-\frac{1}{2+\epsilon}}.
    \]
\end{manualtheorem}
The inequality in Theorem \ref{thm:main rates 2} is best possible in the same sense as the one appearing in Theorem \ref{thm:main rates}.

Another Euclidean quantity that has a prominent role in the literature is $\lvert f^n(z)-\tau\rvert$ (see, for example, \cite[Chapter 16]{BCDM-Book}). Simple arguments in hyperbolic geometry allow us to obtain the next corollary of Theorem \ref{thm:main rates}.

\begin{cor}\label{cor:main rates cor}
    Let $f\colon\D\to\D$ be a non-elliptic map with Denjoy--Wolff point $\tau\in\partial\D$ and whose Koenigs domain is not the whole complex plane. For any $z\in\D$, we have that
    \[
    \limsup\limits_n\frac{\log|f^n(z)-\tau|}{\log n}\le -\frac{1}{4}.
    \]
    If, in addition, $\{f^n(z)\}$ converges to $\tau$ non-tangentially for some $z\in\D$, then $\tfrac{1}{4}$ can be replaced by $\tfrac{1}{2}$. 
\end{cor}

In the special case where the boundary of the Koenigs domain of $f$ has positive logarithmic capacity, we prove a sharper estimate for the Euclidean rate $\lvert f^n(z)-\tau\rvert$ that will be stated in Theorem \ref{thm:rate-zerostep-nonpolar}. The arguments of this result combine the semigroup-fication with estimates for the harmonic measure, and are inspired by a result of Betsakos \cite[Theorem 1]{Betsakos-Rate-Par}.

In order to prove Theorem \ref{thm:main rates}, we employ the semigroup-fication technique to study the geometry of domains $\Omega\subsetneq\C$ satisfying $\Omega+1\subset\Omega$, which are not necessarily asymptotically starlike at infinity. Such a domain always carries a hyperbolic distance $d_\Omega$, for which we prove the following estimate. 

\begin{prop}\label{prop:estimate on Koenigs domains}
    Let $\Omega\subsetneq\C$ be a domain satisfying $\Omega+1\subset\Omega$. For any $z\in\Omega$, we have that
    \begin{equation}\label{eq:estimate on Koenigs domains eq}
        \liminf_n \frac{d_\Omega(z,z+n)}{\log n}\geq \frac{1}{4}.
    \end{equation}
\end{prop}

As a simple example of the domains described by Proposition \ref{prop:estimate on Koenigs domains}, one can think of $\Omega_\N\vcentcolon=\C\setminus \{-n\colon n\in\N\}$. Of course, a generic domain of this type can be vastly more complicated and thus its hyperbolic geometry is particularly difficult to handle directly. This is evident by the lack of estimates similar to \eqref{eq:estimate on Koenigs domains eq} in the literature; even for (seemingly) simple cases  such as $\Omega_\N$.

For $\Omega_\N$ in particular, we prove that the limit inferior in \eqref{eq:estimate on Koenigs domains eq} is in fact a limit (see Proposition \ref{prop:extremal domain} and \eqref{eq:discrete Koebe rate}). As such, the estimate in Proposition \ref{prop:estimate on Koenigs domains} is sharp. 

We end the Introduction with an application of Theorem \ref{thm:main rates} to operator theory. Each holomorphic map $f\colon\D\to\D$ induces a \emph{composition operator} $C_f\colon X\to X$ with $C_f(g)=g\circ f$, where $X$ is either the Hardy space $H^p$, for $p\geq1$, or the Bergman space $A_\alpha^p$, for $p\geq1$ and $\alpha>-1$, in the unit disc. Littlewood's Subordination Principle tells us that such a composition operator is bounded. The theory of composition operators is often intertwined with holomorphic dynamics \cites{Betsakos-Hardy,BMS,BGGY,CZZ}.

 Observe that the operator $C_f^{\, n}\vcentcolon=C_{f^n}$ is also well-defined and bounded, and write $\lVert C_f^{\, n}\rVert_{H^p}$ and $\lVert C_f^{\, n}\rVert_{A^p_\alpha}$ for the norms of $C_f^{\, n}$ in $H^p$ and $A_\alpha^p$, respectively. A result of Arosio and Bracci \cite[Proposition 5.8]{Arosio-Bracci} shows that the limits 
\[
\ell_p\vcentcolon=\lim_{n\to+\infty}\frac{\log\lVert C_f^{\, n}\rVert_{H^p}}{n}, \quad \ell_{p,\alpha}\vcentcolon=\lim_{n\to+\infty}\frac{\log\lVert C_f^{\, n}\rVert_{A^p_\alpha}}{ n},
\]
exist for any non-elliptic $f\colon\D\to\D$. The existence of $\ell_p$ and $\ell_{p,\alpha}$ also follows from standard operator-theoretic arguments, since these quantities are the logarithms of the spectral radii of the operator $C_f$ in $H^p$ and $A^p_\alpha$, respectively (see, for example, \cite[Theorem 3.9]{CMC} for the case of $\ell_p$). 

In particular, if $f$ is hyperbolic then $\ell_{p,\alpha}=(2+\alpha)\ell_p>0$, while if $f$ is parabolic $\ell_{p,\alpha}=\ell_p=0$ (see Corollary \ref{cor:growth bounds}). Using Theorem \ref{thm:main rates} we show that for a parabolic $f$, a more precise estimate for the growth rates of $\lVert C_f^{\, n}\rVert_{H^p}$ and $\lVert C_f^{\, n}\rVert_{A^p_\alpha}$ is attainable.

\begin{cor}\label{cor:operators}
    Let $f\colon\D\to\D$ be a non-elliptic map whose Koenigs domain is not the whole complex plane. For all $p\geq1$ and $\alpha>-1$, we have that
    \begin{enumerate}[label=\rm (\alph*)]
    \item $\displaystyle  \liminf\limits_n\frac{\log\lVert C_f^{\, n}\rVert_{H^p}}{\log n}\geq \frac{1}{2p}$;  and
    \item $\displaystyle \liminf\limits_n\frac{\log\lVert C_f^{\, n}\rVert_{A^p_\alpha}}{\log n}\geq \frac{2+\alpha}{2p}$. 
    \end{enumerate}
\end{cor}
The inequalities in Corollary \ref{cor:operators} are sharp, due to the sharpness of Theorem \ref{thm:main rates}.

\subsection*{Structure of the article.} In Section \ref{sect:prelim} we review two concepts from complex analysis relevant to our work. Additional information on holomorphic iteration, as well as the basic concepts from the theory of one-parameter semigroups can be found in Section \ref{sect:dynamics}. Section \ref{sect: hyperbolic geometry} contains the basics of hyperbolic geometry.

Our extension of the Bracci--Roth semigroup-fication is spread across Sections \ref{sect:internal tangency} and \ref{section:fundamental domains and semigroupfication}. In particular, in Section \ref{sect:internal tangency} we develop a theory for two simply connected domains $D_1\subset D_2$, whose boundaries are similar close to some prime end $\zeta$ of $D_2$. We show that in such a scenario, the hyperbolic distances $d_{D_1}$ and $d_{D_2}$ are Lipschitz equivalent close to $\zeta$. These results might be of independent interest. The constructions involved in Theorems \ref{thm:main A} and \ref{thm:main B} are realised in Section \ref{section:fundamental domains and semigroupfication}.

Section \ref{sect:slope} contains the proofs of Theorem \ref{thm:main slope} and its corollary, Corollary \ref{coro:main slope coro}, while our result on the rate of convergence, Theorem \ref{thm:main rates}, and its consequences are proved in Section \ref{section:convergence rates}. Section \ref{section:convergence rates} also includes lower bounds on the Euclidean rate of convergence, which do not appear in the literature, but can be easily derived by known results. Finally, Section \ref{sect:composition operator} contains a proof of Corollary \ref{cor:operators}.

%%%%%%%%%%%%%%%%%%%%%%%%%%%
\section{Preliminaries}\label{sect:prelim}
%%%%%%%%%%%%%%%%%%%%%%%%%%%

\subsection{Boundaries of simply connected domains}\label{sect:Caratheodory}
Our analysis often requires us to discuss the manner in which sequences or curves approach the boundary of a simply connected domain. For this we turn to the theory of prime ends (see \cite[Chapter 4]{BCDM-Book} or \cite[Chapter 2]{Pommerenke} for a complete presentation).

Consider the extended complex plane $\widehat{\C}\vcentcolon=\C\cup\{\infty\}$ equipped with the spherical metric. Let $D\subsetneq\C$ be a simply connected domain and $\gamma\vcentcolon[0,1]\to\widehat{\C}$ a Jordan arc. The trace $C\vcentcolon=\gamma([0,1])$ is called a \textit{cross cut} of $D$ if $\gamma((0,1))\subset D$ and $\gamma(0), \gamma(1)\in\partial_\infty D$, where $\partial_\infty D$ is the boundary of $D$ in $\widehat{\C}$. When $C$ is a cross cut of $D$, the open set $D\setminus C$ consists of two open connected components $A$ and $B$ satisfying $\partial A\cap D=\partial B\cap D=C\cap D$. A \textit{null-chain} of $D$ is a sequence of cross cuts $\{C_n\}$ that satisfies the following three conditions:
\begin{enumerate}
    \item[(i)] $C_n\cap C_m=\emptyset$, for all $n,m\in\N$, $n\ne m$;
    \item[(ii)] for each $n\geq2$, the sets $C_1\cap D$ and $C_{n+1}\cap D$ lie in different connected components of $D\setminus C_n$;
    \item[(iii)] the spherical diameter of $C_n$ converges to $0$, as $n\to+\infty$.
\end{enumerate}
Given a null-chain $\{C_n\}$ and $n\geq2$, the \textit{interior part} of $C_n$ is defined as the connected component of $D\setminus C_n$ that does not contain $C_1\cap D$. We use $V_n$ to denote the interior part of $C_n$. Two null-chains $\{C_n\}$ and $\{C_n'\}$ are said to be \textit{equivalent} if for every $n\geq2$ there exists $m\in\N$ so that
\begin{equation*}
    V_m'\subset V_n \quad \textup{and} \quad V_m\subset V_n',
\end{equation*}
where $V_n'$ is the interior part of $C_n'$. This is indeed an equivalence relation on the null-chains of $D$. An equivalence class is called a \emph{prime end of $D$}, and the set of all equivalence classes is denoted by $\partial_C D$. The \emph{impression} of a prime end $\zeta$ of $D$, represented by a null-chain $\{C_n\}$, is the non-empty set
\begin{equation}\label{eq:impression of prime end}
    I(\zeta)\vcentcolon=\bigcap\limits_{n=2}^{+\infty}\overline{V_n},
\end{equation}
where the closures are taken in $\widehat{\C}$. It is easy to see that the impression is independent of the choice of the null-chain. 

The prime ends of a simply connected domain $D\subsetneq\C$ induce a topology on $D\cup \partial_CD$, called the Carath\'eodory topology, that agrees with the usual topology of $D$. This topology takes its name from a theorem of Carath\'eodory which shows that any Riemann map $f\colon \D\to D$ can be extended to a homeomorphism $f\colon \D\cup\partial\D\to D\cup \partial_CD$. As a slight abuse of notation we use the same symbol for the Riemann map and its Carath\'eodory extension. 

For the rest of this subsection, let $D\subsetneq\C$ be a simply connected domain and $f\colon \D\cup\partial\D\to D\cup \partial_CD$ the Carath\'eodory extension of a Riemann map. Given $\sigma\in\partial\D$ and $R>1$, the set
\begin{equation}\label{eq:Stolz}
    S(\sigma,R)\vcentcolon=\left\{z\in\D\vcentcolon\dfrac{\lvert\sigma-z\rvert}{1-\lvert z\rvert}<R\right\}
\end{equation}
is called a \textit{Stolz angle} of the unit disk at $\sigma$. A sequence $\{z_n\}\subset\D$ converging to $\sigma\in\partial\D$ is said to converge to \emph{non-tangentially} if there exists $R>1$ such that $\{z_n\}\subset S(\sigma,R)$. Throughout the text we follow the terminology described below. 

\begin{defin}\label{defin:convergence in caratheodory topology}
Let $\zeta\in\partial_CD$ and suppose that $\sigma\in\partial\D$ is the unique point with $f(\sigma)=\zeta$. Consider a sequence $\{w_n\}\subset D$ and a curve $\gamma\colon[0,+\infty)\to D$.
\begin{enumerate}[label=(\roman*)]
    \item We write that $\lim_{n\to+\infty}w_n=\zeta$ \textit{in the Carath\'{e}odory topology of} $D$ provided that $\lim_{n\to+\infty}f^{-1}(w_n)=\sigma$.
    \item We say that $\{w_n\}$ converges to $\zeta$ \textit{non-tangentially in $D$} if and only if $\{f^{-1}(w_n)\}$ converges to $\sigma$ non-tangentially. 
    \item We say that $\gamma$ \textit{lands at} $\zeta$ if $\lim_{t\to+\infty}f^{-1}(\gamma(t))=\sigma$ in the Euclidean topology of $\D$. In addition, $\gamma$ lands at $\zeta$ \emph{non-tangentially} if the curve $f^{-1}\circ\gamma$ is contained in a Stolz angle at $\sigma$. 
\end{enumerate}
\end{defin}

The next definition extends the notion of the slope presented in the Introduction.

\begin{defin}\label{def:slope}
 Fix a prime end $\zeta\in\partial_C D$ and denote by $\sigma$ the unique point of $\partial\D$ with $f(\sigma)=\zeta$. Let $\{z_n\}\subset D$ be a sequence converging to $\zeta$ in the Carath\'{e}odory topology of $D$. Its \textit{slope in $D$}, is defined as $\textup{Slope}_D(z_n)\vcentcolon =\textup{Slope}_\D(f^{-1}(z_n))$. The slope $\textup{Slope}_D(\gamma)$ of a curve $\gamma:[0,+\infty)\to D$ landing at $\zeta$ is defined similarly.
\end{defin}

Note that, by definition, the slope is a conformally invariant quantity. Furthermore, the slope is always a non-empty subset of $[-\tfrac{\pi}{2},\tfrac{\pi}{2}]$. Particularly for curves, $\textup{Slope}_D(\gamma)$ is a continuum.

It is easy to see that a sequence $\{z_n\}\subset D$ converges to $\zeta\in\partial_CD$ non-tangentially if and only if $\textup{Slope}_D(z_n)\subset(-\tfrac{\pi}{2},\tfrac{\pi}{2})$ and similarly for a curve of $D$ landing at $\zeta$. On the other hand, we say that $\{z_n\}$ converges to $\zeta$ \textit{tangentially} if $\textup{Slope}_D(z_n)\subset\{-\tfrac{\pi}{2},\tfrac{\pi}{2}\}$. If $\gamma$ is a curve of $D$ that lands at $\zeta$, we say that it lands \emph{tangentially} if $\textup{Slope}_D(\gamma)=\{-\tfrac{\pi}{2}\}$ or $\textup{Slope}_D(\gamma)=\{\tfrac{\pi}{2}\}$ (the connectedness of $\textup{Slope}_D(\gamma)$ implies that it cannot contain both $-\tfrac{\pi}{2}$ and $\tfrac{\pi}{2}$). Let us emphasise that the absence of non-tangential convergence is different from tangential convergence.

\subsection{Harmonic measure}
One of our results on the rate of convergence requires techniques involving the harmonic measure. All the information presented in this subsection can be found in \cites{Beliaev}, or \cite[Chapter 7]{BCDM-Book}, \cite[Section 4.4]{Pommerenke} for the case of simply connected domains. 

Let $D\subset\C$ be a domain whose Euclidean boundary $\partial D$ is non-polar; i.e. has positive logarithmic capacity. Let $E$ be a Borel subset of $\partial D$. Then, the \textit{harmonic measure} of $E$ with respect to $D$ is exactly the solution of the generalized Dirichlet problem
$$\begin{cases}
	\Delta u=0 \quad \text{in }D,\\
	u=\chi_E \quad\text{on }\partial D.
\end{cases}$$
For $z\in D$ we will use $\omega(z,E,D)$ to denote this solution. By definition, $\omega(\cdot,E,D)$ is a harmonic function on $D$ for every choice of Borel set $E\subset\partial D$, while $\omega(z,\cdot,D)$ is a Borel probability measure on $\partial D$, for each $z\in D$.

An important aspect of the harmonic measure is that it satisfies a subordination principle. To describe this, consider two domains $D_1,D_2$ with non-polar boundaries, and Borel sets $E_1\subset\partial D_1$ and $E_2\subset\partial D_2$. Let $f\vcentcolon D_1\to D_2$ be a holomorphic map that extends continuously (in Euclidean terms) to $E_1$, with $f(E_1)\subset E_2$. Then
\begin{equation}\label{eq:harmonic measure subordination}
	\omega(z,E_1,D_1)\le\omega(f(z),E_2,D_2), \quad \text{for all}\ z\in D_1,
\end{equation}
with equality if and only if $f$ is a homeomorphism between $D_1\cup E_1$ and $D_2\cup E_2$.

The subordination principle yields a domain monotonicity property. That is, given $D_1\subset D_2$ with non-polar boundaries and a Borel set $E\subset\partial D_1\cap\partial D_2$, we have
\begin{equation}\label{eq:harmonic measure monotonicity}
	\omega(z,E,D_1)\le\omega(z,E,D_2), \quad\text{for all}\ z\in D_1.
\end{equation}

Particularly for the case of the unit disc, whenever $E\subset\partial\D$ is an arc, simple calculations lead to the handy formula
\begin{equation}\label{eq:harmonic measure arc}
    \omega(0,E,\D)=\frac{1}{\pi}\arcsin\left(\frac{\textup{diam}[E]}{2}\right).
\end{equation}
In this setting we also have the following result.

\begin{theorem}[\cites{FRW,Solynin}]\label{thm:solynin}
	Let $E\subset\overline{\D}\setminus\{0\}$ be a Jordan arc and let $d:=\textup{diam}[E]$. Denote by $D$ the connected component of $\D\setminus E$ that contains $0$. Let $E_d$ be an arc on $\partial\D$ satisfying $\textup{diam}[E_d]=d$ (in the extremal case when $d=2$, we take $E_d$ to be a half-circle). Then
	\begin{equation}\label{eq:solynin}
		\omega(0,E,D)\ge\omega(0,E_d,\D).
	\end{equation}
\end{theorem}

%%%%%%%%%%%%%%%%%%%%%%%%%%%%%%%%%
\section{Holomorphic dynamics}\label{sect:dynamics}
%%%%%%%%%%%%%%%%%%%%%%%%%%%%%%%%%

\subsection{Iteration in the unit disc}
We now present certain supplementary material from the theory of holomorphic iteration. For further details see \cite{Abate}. 

Recall that we denote by $\tau\in\partial\D$ the Denjoy--Wolff point of a non-elliptic, holomorphic map $f\colon \D\to\D$. The Julia--Carath\'eodory Theorem implies that \emph{the angular derivative} $f'(\tau)\vcentcolon=\angle\lim_{z\to\tau}f'(z)$ of $f$ at $\tau$ exists and satisfies $f'(\tau)\in(0,1]$. This gives a first classification of non-elliptic maps; that is, $f$ is hyperbolic if $f'(\tau)<1$ and it is parabolic otherwise.

Another important aspect of the behaviour of $f$ close to $\tau$ is given by Julia's Lemma, which states that
\begin{equation}\label{eq:Julia's lemma}
   \frac{\lvert \tau-f(z)\rvert^2}{1-\lvert f(z)\rvert^2}\leq f'(\tau) \frac{\lvert \tau-z\rvert^2}{1-\lvert z\rvert^2}, \quad \text{for all}\ z\in\D.
\end{equation}

As mentioned in the Introduction, a Koenigs domain $\Omega$ and a Koenigs function $h$ of a self-map $f$ are unique up to translation, but it is important to note that this essential uniqueness relies on the fact that we require $\Omega$ to be asymptotically starlike at infinity (see \cite{Arosio-Bracci}, or \cite[Section 3.5]{Abate}). When looking at pairs $h,\Omega$ that satisfy
\begin{equation}\label{eq:Koenigs}
 h(\D)=\Omega\quad \text{and} \quad h(f(z))=h(z)+1, \quad \text{for all }z\in\D,   
\end{equation}
but with no further assumptions on $\Omega$, this uniqueness may fail (see \cite{CDP2}).

When $f$ is univalent, $\Omega$ is simply connected and $h$ is simply a Riemann map. As of yet, it is not known whether $\Omega$ can be the whole complex plane.

For the rest of this article, we use the term \emph{positive-parabolic} to describe an $f$ that is parabolic of positive hyperbolic step  and the term \emph{zero-parabolic} for an $f$ that is parabolic of zero hyperbolic step.

The type of a non-elliptic map has important implications on the slope with which its orbits approach the Denjoy--Wolff point. For a hyperbolic $f\colon\D\to\D$, Wolff \cite{Wolff} showed that for each $z\in\D$ there exists some $\theta_z\in(-\tfrac{\pi}{2},\tfrac{\pi}{2})$ so that $\textup{Slope}_\D(f^n(z))=\{\theta_z\}$. This implies that each orbit $\{f^n(z)\}$ converges to $\tau$ non-tangentially. Next, for a positive-parabolic $f$  we have that either $\textup{Slope}_\D(f^n(z))=\{-\tfrac{\pi}{2}\}$ for all $z\in\D$ or $\textup{Slope}_\D(f^n(z))=\{\tfrac{\pi}{2}\}$ for all $z\in\D$; see \cite[Remark 2.3]{CCZRP} and \cite{Pommerenke - Iteration}. In any case, all orbits of positive-parabolic maps converge tangentially. For zero-parabolic maps, the situation is far more chaotic. In \cite{CCZRP} the authors show that the slope of $\{f^n(z)\}$ is independent of the choice of $z\in\D$; that is $\textup{Slope}_\D(f^n(z_1))=\textup{Slope}_\D(f^n(z_2))$, for all $z_1,z_2\in\D$. They also prove that given any compact, connected set $\Theta\subset[-\frac{\pi}{2},\frac{\pi}{2}]$, there exists a zero-parabolic map $f:\D\to\D$ with $\textup{Slope}_\D(f^n(z))=\Theta$, for any $z\in\D$.

Thus, either all orbits of $f$ converge non-tangentially to the Denjoy--Wolff point or none does. So, for simplicity, we say that $\{f^n\}$ converges to $\tau$ non-tangentially whenever some orbit $\{f^n(z)\}$ converges non-tangentially.

\subsection{Semigroups of holomorphic functions}\label{sect:semigroups}
The theory of one-parameter semigroups of holomorphic functions has flourished in the past thirty years, with many of its advances being influential in fields such as geometric function theory, operator theory and the theory of conformal invariants, to name a few. For a complete presentation of this topic containing most recent results, we refer to \cite{BCDM-Book}.

Even though semigroups are typically studied in the context of the unit disc, the majority of their theory remains valid in any simply connected domain $D\subsetneq\C$. As such, we say that a family $(\phi_t)$, for $t\geq0$, of holomorphic functions $\phi_t\vcentcolon D\to D$, is a \emph{semigroup in $D$} if 
\begin{enumerate}[label=\rm(\roman*)]
    \item $\phi_0=\mathrm{Id}_D$;
    \item $\phi_{t+s}=\phi_t\circ \phi_s$, for all $t,s\geq 0$; and
    \item $\lim_{t\to0}\phi_t(z)=z$, for all $z\in D$.
\end{enumerate}
An important consequence of the definition is that $\phi_t$ is univalent for each $t\ge 0$ (see \cite[Theorem 8.1.17]{BCDM-Book}).

If $(\phi_t)$ is a semigroup in $D$, the theory of prime ends and the continuous version of the Denjoy--Wolff Theorem \cite[Theorem 8.3.1]{BCDM-Book} imply that there exists a unique $\tau\in D\cup\partial_C D$ such that $\lim_{t\to+\infty}\phi_t(z)=\tau$, in the Carath\'{e}odory topology of $D$, for all $z\in D$. If $\tau\in\partial_CD$, we say that $(\phi_t)$ is \textit{non-elliptic}, and $\tau\in\partial_C D$ is called the \textit{Denjoy--Wolff prime end} of $(\phi_t)$. When the Euclidean boundary of $D$ is ``simple" (in the case of a disc, for instance), the Denjoy--Wolff prime end is merely a point, which we call \emph{the Denjoy--Wolff point of $(\phi_t)$}. An analysis on the difference between ``attracting" prime ends and points can be found in \cite{Bracci-Benini}.

The linearisation through a Koenigs function we described in the case of discrete iteration extends to semigroups. That is, for a non-elliptic semigroup $(\phi_t)$ in $D$, there exists a \emph{Koenigs domain} $\Omega\subsetneq\C$ that is starlike at infinity and a \emph{Koenigs function} $h: D\to \Omega$, with $h(D)=\Omega$, such that
\begin{equation}\label{eq:Koenigs semigroup}
    h(\phi_t(z))=h(z)+t, \quad\textup{for all }z\in D \textup{ and all }t\ge0.
\end{equation}
The pair $h,\Omega$ is uniquely determined up to translation. A Koenigs domain of a semigroup is always simply connected and a Koenigs function is always univalent.

There are three mutually exclusive possibilities for the simply connected domain $\bigcup_{t\ge0}(\Omega-t)$: a horizontal strip, a horizontal half-plane, and all of $\C$. Just like before, we say that the semigroup is \textit{hyperbolic}, \textit{positive-parabolic} and \textit{zero-parabolic}, respectively in these three cases.

For a non-elliptic semigroup $(\phi_t)$ in a simply connected domain $D$, we say that the curve $\eta_z\colon[0,+\infty)\to D$ with $\eta_z(t)=\phi_t(z)$, for some $z\in D$, is a \emph{trajectory} of $(\phi_t)$. We often denote the trajectory $\eta_z$ by $(\phi_t(z))$, for simplicity. By the continuous version of the Denjoy--Wolff Theorem we mentioned, each trajectory lands at the Denjoy--Wolff prime end $\tau\in\partial_CD$. The slope $\mathrm{Slope}_D(\eta_z)$ of a trajectory is in a one-to-one analogy with the discrete setting. That is, if $(\phi_t)$ is hyperbolic, then $\textup{Slope}_D(\eta_z)=\{\theta_z\}$ with $\theta_z\in(-\frac{\pi}{2},\frac{\pi}{2})$ depending on $z\in D$. When $(\phi_t)$ is positive-parabolic, either $\textup{Slope}_D(\eta_z)=\{-\frac{\pi}{2}\}$ for all $z\in D$ or $\textup{Slope}_D(\eta_z)=\{\frac{\pi}{2}\}$ for all $z\in D$. Finally, when $(\phi_t)$ is zero-parabolic all trajectories have the same slope, which can be any continuum in $[-\frac{\pi}{2},\frac{\pi}{2}]$; for more information see \cites{BCDM-Book, CDM, Kelgiannis}.

The rate of convergence of semigroups has been the topic of extensive research over the past twenty years (see, for example, \cites{Betsakos-Rate-Hyp,Betsakos-Rate-Par,BCDM-Rates,BCZZ,Bracci-Speeds}). The result most relevant to our analysis is the following, taken from \cite[Theorem 16.3.3]{BCDM-Book}.

\begin{theorem}\label{thm:semigroup rates}
    Let $(\phi_t)$ be a non-elliptic semigroup in $\D$ with Denjoy--Wolff point $\tau\in\partial\D$. For every $z\in\D$ there exists a positive constant $c\vcentcolon=c(z)$ so that 
    \[
    \lvert \phi_t(z)-\tau\rvert \leq \frac{c}{\sqrt{t}},\quad \text{for all}\ t>0.
    \]
\end{theorem}

%%%%%%%%%%%%%%%%%%%%%%%%%%%%%%%%%%%%%%%%%%%%%%%%%%%%%%%%%%%%%%%
\section{Hyperbolic geometry}\label{sect: hyperbolic geometry}
%%%%%%%%%%%%%%%%%%%%%%%%%%%%%%%%%%%%%%%%%%%%%%%%%%%%%%%%%%%%%%%

We now present the main aspects of hyperbolic geometry required for our techniques. For more information we refer to the books \cite[Chapter 1]{Abate},\cite[Chapter 5]{BCDM-Book}.

We are interested in domains $D\subset\C$ for which $\C\setminus D$ contains at least two points. For such a domain, called a \emph{hyperbolic domain}, there exists a universal covering $\pi\colon \D\to D$, i.e. a local biholomorphism that has the path lifting property, which is unique up to pre-composition with a M\"obius automorphism of the unit disc. When $D$ is a simply connected domain, other than the complex plane, $\pi$ is a Riemann map.

There is a unique metric $\lambda_D(z)\lvert dz\rvert$ in $D$, called the \emph{hyperbolic metric of $D$}, that is independent of the choice of the universal covering and satisfies
\begin{equation}\label{eq:hyperbolic metric in domain}
\lambda_D(\pi(\tilde{z}))\lvert \pi'(\tilde{z})\rvert \rvert d\tilde{z}\rvert = \lambda_\D(\tilde{z}) \rvert d\tilde{z} \lvert , \quad \text{for all} \ \tilde{z}\in\D,
\end{equation}
where the quantity 
\begin{equation}\label{eq:hyperbolic metric in D}
    \lambda_\D(\tilde{z}) \lvert d\tilde{z}\rvert =\frac{\rvert d\tilde{z} \lvert}{1-|\tilde{z}|^2},
\end{equation}
is the \emph{hyperbolic metric of the unit disc}. 

The \emph{hyperbolic length} of a piecewise $C^1$-smooth curve $\gamma:I\to D$, defined on an interval $I\subseteq\R$, is
\begin{equation}\label{eq:hyperbolic length in a domain}
	\ell_D(\gamma)=\int\limits_{\gamma}\lambda_D(z)|dz|=\int\limits_I\lambda_D(\gamma(t))|\gamma'(t)|dt.
\end{equation}  
In addition, for $t_1\le t_2$ in $I$, we define
\begin{equation}\label{eq:hyperbolic length between two points}
	\ell_D(\gamma;[t_1,t_2])=\int\limits_{t_1}^{t_2}\lambda_D(\gamma(t))|\gamma'(t)|dt.
\end{equation}

Equation \eqref{eq:hyperbolic metric in domain} implies that if $\gamma\colon I\to D$ is a piecewise $C^1$-smooth curve and $\tilde{\gamma}\colon I\to\D$ a \emph{lift of $\gamma$}; i.e. a curve satisfying $\pi\circ\tilde{\gamma}=\gamma$, then
\begin{equation}\label{eq:universal cover is a local isometry}
    \ell_D(\gamma)=\int_I\lambda_D(\gamma(t))\lvert \gamma'(t)\rvert dt=\int_I\lambda_\D(\tilde{\gamma}(t)) \lvert \tilde{\gamma}'(t)\rvert dt = \ell_\D(\tilde{\gamma}).
\end{equation}

The \textit{hyperbolic distance} between $z,w\in D$ is defined as 
\begin{equation}\label{eq:hyperbolic distance in domain}
	d_D(z,w)=\inf\limits_{\gamma}\ell_D(\gamma),
\end{equation}
where the infimum is taken over all piecewise $C^1$-smooth curves $\gamma$ in $D$ joining $z$ and $w$. Particularly for the case of $\D$ the hyperbolic distance can be computed explicitly and has the following closed-form formula:
\begin{equation}\label{eq:hyperbolic distance unit disk}
	d_\D(z,w)=\frac{1}{2}\log\frac{|1-\bar{w}z|+|z-w|}{|1-\bar{w}z|-|z-w|}, \quad\quad z,w\in\D.
\end{equation}

The universal covering also satisfies 
\begin{equation}\label{eq:universal cover is a contraction}
    d_D(\pi(\tilde{z}),\pi(\tilde{w}))\leq d_\D(\tilde{z},\tilde{w}), \quad \text{for all}\ \tilde{z},\tilde{w}\in\D,
\end{equation}
When $\pi$ is a Riemann map we have equality in \eqref{eq:universal cover is a contraction} for all $\tilde{z},\tilde{w}\in\D$.

Let us now introduce some notation in a  hyperbolic domain $D$. For a point $w\in D$ and $R>0$, we use $D_D(w,R)$ to denote the hyperbolic disk of $D$ centred at $w$ and of radius $R$; that is
\begin{equation}\label{eq:hyperbolic disk}
	D_D(w,R)=\{z\in D:d_D(z,w)<R\}.
\end{equation}
Also, if $\gamma\colon I\to D$ is a curve defined on some interval $I\subseteq \R$ and $z\in D$, for convenience we write
\[
d_D(z,\gamma)\vcentcolon=\inf\{d_D(z,\gamma(t))\colon t\in I\}.
\]

One of the main advantages of working with the hyperbolic metric is an extension of the Schwarz--Pick Lemma stating that if $f\colon D_1\to D_2$ is a holomorphic map between hyperbolic domains $D_1,D_2$, then 
\begin{equation}\label{eq:hyperbolic contraction}
	d_{D_2}(f(z),f(w))\le d_{D_1}(z,w), \quad z,w\in D_1.
\end{equation}
This gives rise to the \textit{domain monotonicity property} of the hyperbolic distance, where if $D_1\subset D_2$ then
\begin{equation}\label{eq:hyperbolic domain monotonicity}
	d_{D_2}(z,w)\le d_{D_1}(z,w), \quad \text{for all } z,w\in D_1.
\end{equation}

We now describe the geodesics of hyperbolic domains. Focusing on $\D$ momentarily, a curve $\gamma\colon I\to\D$ defined on an interval $I\subseteq \R$ is called a \emph{(hyperbolic) geodesic of $\D$} if for any $t_1\leq t_2$ in $I$ we have that 
\[
\ell_\D(\gamma;[t_1,t_2])=d_\D(\gamma(t_1),\gamma(t_2)). 
\]
Since we will not be considering any type of geodesic other than a hyperbolic geodesic, in most cases we will omit the term ``hyperbolic". Furthermore, we sometimes also refer to the trace of $\gamma$ when using the term geodesic. The geodesics of $\D$ are parts of circles or straight lines that are perpendicular to the unit circle $\partial \D$. Hence, every two distinct points $z,w\in\D$ can be joined by a unique geodesic of $\D$.

In a general hyperbolic domain $D$ the geodesics are exactly the curves that lift to geodesics of the unit disc. That is, a curve $\gamma\colon I\to D$, for some interval $I\subseteq\R$, is a \emph{(hyperbolic) geodesic} of a hyperbolic domain $D$, if there exists a geodesic $\tilde{\gamma}\colon I\to\D$ of $\D$ so that $\pi\circ\tilde{\gamma}=\gamma$.

When $D$ is a simply connected hyperbolic domain, the conformal invariance of the hyperbolic distance implies that for every $z,w$ there exists a unique geodesic $\gamma\colon [0,1]\to D$ of $D$ joining $z$ and $w$, which also satisfies $d_D(z,w)=\ell_D(\gamma)$. Also, if $\gamma\colon [0,+\infty)\to D$ is a geodesic of $D$ satisfying $\lim_{t\to+\infty}d_D(\gamma(0),\gamma(t))=+\infty$, then $\gamma$ lands at some prime end of $D$.

The situation in multiply connected domains is much more subtle. It turns out that if $D$ is a multiply connected hyperbolic domain and $z,w\in D$ distinct points, there are infinitely many geodesics of $D$ joining $z$ and $w$. Moreover, if $\gamma_1,\gamma_2$ are any two such geodesics of $\Omega$ joining $z$ and $w$, then $\gamma_1$ is not homotopic to $\gamma_2$ in $\Omega$, meaning that we cannot continuously deform $\gamma_1$ to $\gamma_2$ inside the domain $\Omega$ (see \cite[Theorem 7.2.4]{KeenLakic}). It is also known that the geodesics of $D$ might not be minimisers of the hyperbolic distance (see, for example, \cite[Proposition 1.9.30]{Abate}). So, we say that a geodesic $\gamma\colon I\to D$ of $D$ is \emph{minimal} if for any $t_1\leq t_2$ in $I$ we have
\[
\ell_D(\gamma;[t_1,t_2])=d_D(\gamma(t_1),\gamma(t_2)).
\]
Every pair of points $z,w\in D$ can be connected by a minimal geodesic (see \cite[Proposition 1.9.29]{Abate}).

As an example, consider the hyperbolic domain $\Omega_\N\vcentcolon= \C\setminus \{-n\colon n\in \N\}$. This can be thought of as an infinitely connected version of a  slit plane. The next lemma desribes a minimal geodesic of $\Omega_\N$ that will prove important for our analysis of the rates of convergence in Section \ref{section:convergence rates}. 

\begin{lm}\label{lm:minimal geodesic in Omega_N}
    The curve $\gamma\colon (-1,+\infty)\to \Omega_\N$ with $\gamma(t)=t$ is a minimal geodesic of $\Omega_\N$.
\end{lm}

\begin{proof}
    Let $-1<a<b<+\infty$. Consider the open half-plane $H_t=\{-1+e^{it}z\colon\ \mathrm{Re}z>0\}$, for $t\in(-\frac{\pi}{2},\frac{\pi}{2})$. Since $H_t\subset \Omega_\N$ and $a,b\in H_t$, a result of J{\o}rgensen \cite{Jorgensen} tells us that any minimal geodesic of $\Omega_\N$ that joins $a,b$ is contained in $H_t$. Since this holds for any $t\in(-\frac{\pi}{2},\frac{\pi}{2})$, we obtain that any minimal geodesic joining $a,b$ lies in $\bigcap_{t\in(-\frac{\pi}{2},\frac{\pi}{2})}H_t=(-1,+\infty)$, meaning that this geodesic is the real interval $[a,b]$. The result follows from the fact that $a,b$ were arbitrary.
\end{proof}

Determining the geodesics of a hyperbolic domain can be quite a difficult endeavour. So it is often convenient to work with the broader class of quasi-geodesics. These curves were defined in \eqref{eq:quasi-geodesic} for the case of $\D$. This definition can be extended to any hyperbolic domain $D$ by simply replacing the hyperbolic distance $d_\D$ with $d_D$ and $\ell_\D$ with $\ell_D$. If we need to emphasise the constants $A$ and $B$ in the definition, we may use the term \textit{($A,B$)-quasi-geodesic}. It is easy to see that a curve $\gamma\colon [0,+\infty)\to D$ is a quasi-geodesic if and only if $\gamma\lvert_{[T,+\infty)}$ is a quasi-geodesic, for some $T> 0$.

The most important result pertaining quasi-geodesics is the Shadowing Lemma, which for the case of simply connected domains has the following form, found in \cite[Theorem 6.9.8]{BCDM-Book}.

\begin{theorem}[\cite{BCDM-Book}]\label{thm:shadowing lemma}
Assume that $D\subsetneq\C$ is a simply connected domain and that $\eta\colon[0,+\infty)\to D$ is an $(A,B)$-quasi-geodesic of $D$. Then, $\eta$ lands at some $\zeta\in\partial_CD$, and there exist a geodesic $\gamma\colon[0,+\infty)\to D$ of $D$ landing at $\zeta$, and a constant $R>0$ depending only on $A$ and $B$, such that
\[
d_D(\eta(t),\gamma)<R,\quad \text{for all}\ t\in [0,+\infty).
\]
\end{theorem}

We end this section by presenting two results for simply connected domains. The first is known as the \emph{``Distance Lemma"} (see \cite[Theorem 5.3.1]{BCDM-Book}). For this, we use $\delta_D(z):=\textup{dist}(z,\partial D)$ to denote the Euclidean distance of a point $z\in D$ from the boundary $\partial D$ of a domain $D\subset \C$.

\begin{lm}\label{lm:distance lemma}
	Let $D\subsetneq\C$ be a simply connected domain and let $z_1,z_2\in D$. Then
	\begin{equation}\label{eq:distance lemma}
		\frac{1}{4}\log\left(1+\frac{|z_1-z_2|}{\min\{\delta_D(z_1),\delta_D(z_2)\}}\right)\le d_D(z_1,z_2)\le \int\limits_{\gamma}\frac{|d z|}{\delta_D(z)},
	\end{equation}
	where $\gamma$ is any piecewise $C^1$-smooth curve in $D$ joining $z_1$ and $z_2$.
\end{lm}

 Next we have \cite[Lemma 1.8.6]{BCDM-Book}.
 
\begin{lm}\label{lm:convergence to the the boundary in D}
    Let $\{z_n\}$, $\{w_n\}$ be sequences in $\D$ and $\zeta\in\partial\D$. Write
    \[
    C\vcentcolon=\limsup_{n}d_\D(z_n,w_n).
    \]
    \begin{enumerate}[label=\rm(\alph*)]
        \item If $C<+\infty$ and $\{z_n\}$ converges to $\zeta$, then so does $\{w_n\}$. If, in addition, $z_n$ converges to $\zeta$ non-tangentially in $\D$, then the same is true for $\{w_n\}$. 
        \item If $C=0$ and the limit $\lim\limits_{n\to+\infty}\mathrm{arg}\left(1-\overline{\zeta}z_n\right)=\theta\in[-\tfrac{\pi}{2},\tfrac{\pi}{2}]$ exists, then $\lim\limits_{n\to+\infty}\mathrm{arg}\left(1-\overline{\zeta}w_n\right)=\theta$.
    \end{enumerate}
\end{lm}

%%%%%%%%%%%%%%%%%%%%%%%%%%%%%%%%%%%%%%%%%%%%%%%%%%%%%%%%%%%%%%%%%%%%%%%%%%%%%%%%%%%%%%%
\section{Internally tangent simply connected domains}\label{sect:internal tangency}
%%%%%%%%%%%%%%%%%%%%%%%%%%%%%%%%%%%%%%%%%%%%%%%%%%%%%%%%%%%%%%%%%%%%%%%%%%%%%%%%%%%%%%%

Here we develop one of the main tools of our analysis, which roughly shows that when two simply connected domains ``look" very similar close to a prime end, their hyperbolic geometries around the prime end are comparable. 

First, we make the following definition. 

\begin{defin}\label{eq:hyperbolic sector}
Let $D\subsetneq\C$ be a simply connected domain and suppose that $\gamma\colon[0,+\infty)\to D$ is a geodesic of $D$. A \textit{hyperbolic sector around} $\gamma$ \textit{of amplitude} $R>0$ in $D$ is the set
\[
	S_D(\gamma,R)=\{z\in D:d_D(z,\gamma)<R\}.
\]
\end{defin}

The definition of a hyperbolic sector is similar to that of a ``hyperbolic approach region" given in \cite[Definition 2.2.5]{Abate}. Also, in \cite[Lemma 2.2.7 (iii)]{Abate} it is shown that when $D=\D$, hyperbolic sectors are equivalent to the standard Stolz regions defined in \eqref{eq:Stolz}.

In most of our results we also use the notation
\[
S_D\left(\gamma\lvert_{[t_0,+\infty)},R\right)=\{z\in D:d_D(z,\gamma\lvert_{[t_0,+\infty)}))<R\},\quad \text{for}\ t_0\ge0,
\]
to denote the \emph{``tail"} of a hyperbolic sector. Observe that the tail of any hyperbolic sector is also a hyperbolic sector. Moreover, for any $0\leq t_0\leq t_1$ we have the following ``monotonicity"
\[
S_D\left(\gamma\lvert_{[t_1,+\infty)},R\right)\subset S_D\left(\gamma\lvert_{[t_0,+\infty)},R\right).
\]

Hyperbolic sectors were explicitly computed for the case of the right half-plane $\H=\{z\in\C\colon\mathrm{Re}\, z>0\}$ in \cite[Lemma 4.4]{BCDG}, stated below.

\begin{lm}[\cite{BCDG}]\label{lm:sectors in the right half-plane}
    Let $\gamma\colon[0,+\infty)\to\H$ be a geodesic of the right half-plane $\H$, with $\gamma([0,+\infty))\subset \R^+$. Then, for any $R>0$ we have
    \[
    S_\H(\gamma,R)=D_\H(\gamma(0),R)\cup\{re^{i\theta}\colon r>\gamma(0), \ \lvert\theta\rvert<\beta\},
    \]
    where $\beta(R)=\beta\in(0,\pi/2)$ satisfies $d_\H(1,e^{i\beta})=R$.
\end{lm}

Hyperbolic sectors around different geodesics landing at the same prime end $\zeta$ are eventually contained in one another. This result is well-known to experts, but we provide a sketch of the proof for the sake of completeness.

\begin{prop}\label{prop:equivalence of sectors}
Let $D\subsetneq\C$ be a simply connected domain and $\zeta\in\partial_CD$. If $\gamma_1,\gamma_2\colon [0,+\infty)\to D$ are geodesics of $D$ landing at $\zeta$, then for every $R>0$, there exist $R_1,R_2>0$ and $t_0\geq0$ so that
\[
S_D\left(\gamma_2\lvert_{[t_0,+\infty)}, R_1\right)\subset S_D\left(\gamma_1, R\right)\subset S_D\left(\gamma_2, R_2\right)
\]
\end{prop}

\begin{proof}
Conjugating with an appropriate Riemann map allows us to assume that $D$ is the right half-plane $\H$ and $\zeta=\infty$, while $\gamma_1([0,+\infty))=[x_1,+\infty)$ and $\gamma_2([0,+\infty))=\{x+i\colon x\in[x_2,+\infty)\}$, for some $x_1,x_2>0$. Then, by Lemma \ref{lm:sectors in the right half-plane} we have that for any $R>0$
\[
S_\H(\gamma_1,R)=D_\H(x_1,R)\cup\{re^{i\theta}\colon r>x_1, \ \lvert\theta\rvert<\beta\},
\]
for some $\beta(R)=\beta\in(0,\pi/2)$ satisfying $d_\H(1,e^{i\beta})=R$. Since the M\"obius map $z\mapsto z+i$ is a hyperbolic isometry of $\H$, we obtain that
\[
S_\H(\gamma_2,R)=D_\H(x_2+i,R)\cup\{re^{i\theta}+i\colon r>x_2, \ \lvert\theta\rvert<\beta\}.
\]
The result now easily follows from elementary arguments in Euclidean geometry.
\end{proof}

\begin{rem}\label{rem: convexity of sectors}
Although not explicitly stated, the proof of \cite[Lemma 4.6]{BCDG} shows that every hyperbolic sector of a simply connected domain $D$ is a hyperbolically convex set; that is, for every $z,w\in S_D(\gamma,R)$, the geodesic of $D$ joining $z$ and $w$ is contained in $S_D(\gamma,R)$.
\end{rem}

Hyperbolic sectors can be used to characterise the notion of non-tangential convergence given in Definition \ref{defin:convergence in caratheodory topology}, as stated in the following result taken from \cite[Proposition 4.5]{BCDG}.

\begin{prop}[\cite{BCDG}]\label{prop: non-tangential convergence}
Let $D\subsetneq\C$ be a simply connected domain and $\{z_n\}\subset D$ a sequence, such that $z_n\to \zeta\in\partial_CD$ in the Carath\'eodory topology of $D$. The sequence $\{z_n\}$ converges to $\zeta$ non-tangentially if and only if there exist a geodesic $\gamma\colon [0,+\infty) \to D$ of $D$ landing at $\zeta$, and a number $R>0$, such that $\{z_n\}$ is eventually contained in the sector $S_D(\gamma,R)$. 
\end{prop}

The following definition is in some sense an extension of a standard notion in complex analysis, that of an \emph{inner tangent} (see, for example, \cite[Definition 2.4.8]{Abate} and \cite[Definition V. 5.1]{GM}).

\begin{defin}\label{def: internal tangency}
Let $D_1\subset D_2\subsetneq\C$ be two simply connected domains and $\zeta\in\partial_CD_2$. We say that $D_1$ is \emph{internally tangent to $D_2$ at $\zeta$}, if there exists a geodesic $\gamma\colon[0,+\infty)\to D_2$ of $D_2$, landing at $\zeta$, such that for any $R>0$, there exists $t_0\ge0$ so that
\[
S_{D_2}(\gamma\lvert_{[t_0,+\infty)},R)\subset D_1.
\]
\end{defin}

\begin{rem}\label{rem: equiv def of internal tangency}
Note that due to Proposition \ref{prop:equivalence of sectors}, Definition \ref{def: internal tangency} is independent of the choice of the geodesic $\gamma$. Moreover, the notion of internally tangent domains is conformally invariant. Also, by Proposition~\ref{prop: non-tangential convergence}, we can immediately see that $D_1$ is internally tangent to $D_2$ at $\zeta$ if and only if $D_1$ eventually contains any sequence of $D_2$ converging non-tangentially to $\zeta$. 
\end{rem}

One can expect that when $D_1$ is internally tangent to $D_2$ at some $\zeta\in\partial_C D_2$, the boundaries of $D_1$ and $D_2$ look very similar close to $\zeta$. This idea is made precise in the next lemma. 

\begin{lm}\label{lm:intersection of prime ends}
 Let $D_1\subset D_2\subsetneq \C$ be two simply connected domains such that $D_1$ is internally tangent to $D_2$ at $\zeta\in\partial_C D_2$. Then, there exists a unique prime end $\zeta_1\in\partial_CD_1$ of $D_1$ with the following property: If $\{C_n\}$ is a null-chain of $D_1$ representing $\zeta_1$, $\gamma\colon[0,+\infty)\to D_2$ is a geodesic of $D_2$ landing at $\zeta\in\partial_CD_2$ and $R>0$, then for all $n\geq2$ there exists $t_n>0$ so that 
 \begin{equation}\label{eq:intersection of prime ends statement eq}
 S_{D_2}\left(\gamma\lvert_{[t_n,+\infty)},R\right)\subset V_n,
 \end{equation}
 where $V_n$ is the interior part of $C_n$. Moreover, $\zeta_1$ has the same impression as $\zeta$. 
\end{lm}

\begin{proof}
    Using a conformal map we can assume that $D_2=\D$ and $\zeta=1\in\partial\D$. Moreover, due to Proposition \ref{prop:equivalence of sectors} we can choose the geodesic $\gamma\colon[0,+\infty)\to \D$ so that $\gamma([0,+\infty))\subset [0,1)$ and $\lim_{t\to+\infty}\gamma(t)=1$. In this setting $S_\D(\gamma, R)$ can be thought of as a standard Stolz angle given by \eqref{eq:Stolz}.\\
    Since $D_1$ is internally tangent to $\D$ at $1$, there exists some $t_0>0$ so that 
    \begin{equation}\label{eq:intersection of prime ends eq1}
        S_{\D}\left(\gamma\lvert_{[t_0,+\infty)},R\right)\subset D_1.
    \end{equation}
   We are going to construct the desired prime end $\zeta_1\in\partial_CD_1$. Let $C(1,r_n)$ be the Euclidean circle centred at 1 and of radius $r_n>0$, where $\{r_n\}$ is a strictly decreasing sequence converging to 0. Omitting the first few terms, if necessary, we have that 
   \[
   C(1,r_n)\cap S_\D\left(\gamma\lvert_{[t_0,+\infty)},R\right)\neq \emptyset,\quad \text{for all}\ n\in\N.
   \]
   Then, due to \eqref{eq:intersection of prime ends eq1} we have that $C(1,r_n)\cap D_1\neq\emptyset$. Let $C_n$ be the connected component of $C(1,r_n)\cap D_1$ that intersects $S_\D\left(\gamma\lvert_{[t_0,+\infty)},R\right)$. Due to our choice of the sequence $\{r_n\}$, we have that $\{C_n\}$ is a null-chain of $D_1$. We are going to show that the prime end $\zeta_1\in\partial_C D_1$ represented by $\{C_n\}$ has the desired properties.\\
   Fix $n\geq 2$ and let $V_n$ be the interior part of $C_n$. Then,
   \begin{equation}\label{eq:intersection of prime ends eq2}
       S_\D\left(\gamma\lvert_{[t_0,+\infty)},R\right)\cap D(1,r_n)\subset V_n,
   \end{equation}
  where $D(1,r_n)$ is the Euclidean disc bounded by $C(1,r_n)$. We can also choose $t_n\geq t_0$ large enough so that 
  \begin{equation}\label{eq:intersection of prime ends eq3}
      S_\D\left(\gamma\lvert_{[t_n,+\infty)},R\right)\subset D(1,r_n).
  \end{equation}
  Note that since
   \[
   S_\D\left(\gamma\lvert_{[t_n,+\infty)},R\right)\subset S_\D\left(\gamma\lvert_{[t_0,+\infty)},R\right),
   \]
   we can combine \eqref{eq:intersection of prime ends eq2} and \eqref{eq:intersection of prime ends eq3} in order to obtain that
   \[
   S_\D\left(\gamma\lvert_{[t_n,+\infty)},R\right)\subset S_\D\left(\gamma\lvert_{[t_0,+\infty)},R\right)\cap D(1,r_n)\subset V_n,
   \]
   as required for \eqref{eq:intersection of prime ends statement eq}. Furthermore, because $V_n\subset D(1,r_n)$, we immediately get that the impression of $\zeta_1$ is the singleton $\{1\}$.\\   
   Now, if $\{C_n'\}$ is any other null-chain representing $\zeta_1$, then since the interior parts of $\{C_n\}$ and $\{C_n'\}$ are eventually contained in one another we get that \eqref{eq:intersection of prime ends statement eq} will also hold for $V_n'$. Finally, if $\widetilde{\zeta_1}\in\partial_C D_1$ is any prime end, different from $\zeta_1$, represented by some null-chain $\{\widetilde{C_n}\}$, then the interior parts of $C_n$ and $\widetilde{C_n}$ are eventually disjoint, meaning that $\zeta_1$ is the unique prime end of $D_1$ satisfying \eqref{eq:intersection of prime ends statement eq}.
\end{proof}

\begin{rem}
    Whenever the simply connected domain $D_1$ is internally tangent to $D_2$ at $\zeta\in\partial_CD_2$, we are going to say that the prime end $\zeta_1\in\partial_C D_1$ given by Lemma \ref{lm:intersection of prime ends} is the prime end of $D_1$ \emph{associated to $\zeta$}.
\end{rem}

The main result of this section, given below, shows that if $D_1$ is internally tangent to $D_2$ at $\zeta$, then the hyperbolic geometries of $D_1$ and $D_2$ are similar close to $\zeta$. This is inspired by the localization results given in \cite[Section 3]{BCDG}.

\begin{theorem}\label{thm:internal tangency hyperbolic equivalence}
	Let $D_1\subset D_2\subsetneq\C$ be two simply connected domains such that $D_1$ is internally tangent to $D_2$ at $\zeta\in\partial_CD_2$. Fix any geodesic $\gamma:[0,+\infty)\to D_2$ of $D_2$ landing at $\zeta$ and a number $K>1$. Then, for any $R>0$ there exists some $t_1\ge 0$ such that
	\begin{enumerate}
		\item[\textup{(a)}] $S_{D_2}(\gamma|_{[t_1,+\infty)},R)\subset D_1$;
		\item[\textup{(b)}] $\lambda_{D_2}(z)\le \lambda_{D_1}(z)\le K \lambda_{D_2}(z)$, for all $z\in S_{D_2}(\gamma|_{[t_1,+\infty)},R)$; and
		\item[\textup{(c)}] $d_{D_2}(z,w)\le d_{D_1}(z,w)\le K d_{D_2}(z,w)$, for all $z,w\in S_{D_2}(\gamma|_{[t_1,+\infty)},R)$.
	\end{enumerate}
\end{theorem}
\begin{proof}
Let $z\in D_2$ and $c>0$. Take a Riemann map $\phi\colon\mathbb{D}\to D_2$, with $\phi(0)=z$. Then, $\phi\left(D_\D(0,c)\right)=D_{D_2}(z,c)$ due to the conformal invariance of the hyperbolic distance. So, using \eqref{eq:hyperbolic metric in domain}, we get
\begin{equation}\label{eq: internal tang proof 0}
\lambda_{D_{D_2}(z,c)}(z)=\frac{1}{\lvert\phi'(0)\rvert}\lambda_{D_\D(0,c)}(0).
\end{equation}
But, the function $z\mapsto (\tanh c) z$ maps $\D$ conformally onto $D_\D(0,c)$, meaning that 
\begin{equation}\label{eq: internal tang proof 00}
    \lambda_{D_\D(0,c)}(0)=\frac{1}{\tanh c}\lambda_\D(0).
\end{equation}
Combining \eqref{eq: internal tang proof 0} with \eqref{eq: internal tang proof 00} we get that
\begin{equation}\label{eq: internal tang proof}
    \lambda_{D_{D_2}(z,c)}(z)=\frac{1}{\lvert\phi'(0)\rvert}\frac{1}{\tanh c}\lambda_\D(0)=\frac{1}{\tanh c}\lambda_{D_2}(z).
\end{equation}
Now, fix a hyperbolic sector $S_{D_2}(\gamma,R)$ and a number $K>1$ as in the statement of the theorem. We can then choose $c>0$ large enough so that $\frac{1}{\tanh c} <K$. Since $D_1$ is internally tangent to $D_2$ at $\zeta$, there exists $t_0\ge0$ so that $S_{D_2}(\gamma\lvert_{[t_0,+\infty)},R)$ is contained in $D_1$.\\
We claim that there exists $t_1\geq t_0$ so that for any $z\in S_{D_2}(\gamma\lvert_{[t_1,+\infty)},R)$, we have that $D_{D_2}(z,c)\subset D_1$.\\
If this were not the case, then there would exist a sequence $\{t_n\}\subset[t_0,+\infty)$, with $t_n\xrightarrow{n\to+\infty}+\infty$, and a sequence of points $\{z_n\}\subset D_2$ with $z_n\in S_{D_2}(\gamma\lvert_{[t_n,+\infty)},R)$, for every $n\in\N$, so that $D_{D_2}(z_n,c)\cap D_1^c\neq\emptyset$. Note that $\{z_n\}$ converges non-tangentially to $\zeta$ in $D_2$. Choose a sequence  $w_n\in D_{D_2}(z_n,c)\cap D_1^c$, for $n\in\N$; that is $d_{D_2}(w_n,z_n)<c$ for all $n\in\N$. According to Lemma \ref{lm:convergence to the the boundary in D} (a), this means that $\{w_n\}$ also converges to $\zeta$ non-tangentially in $D_2$. But, the fact that $\{w_n\}$ is not contained in $D_1$ contradicts the equivalent definition of internally tangent domains given in Remark~\ref{rem: equiv def of internal tangency}.\\
With our claim proved, notice that since $t_1\geq t_0$, we have that $S_{D_2}(\gamma\lvert_{[t_1,+\infty)},R)$ is contained in $ S_{D_2}(\gamma\lvert_{[t_0,+\infty)},R)\subset D_1$, which yields (a). Let us now consider a point $z\in S_{D_2}(\gamma\lvert_{[t_1,+\infty)},R)$. By the domain monotonicity of the hyperbolic metric \eqref{eq:hyperbolic domain monotonicity}, our claim, and \eqref{eq: internal tang proof}, we have that
\[
\lambda_{D_2}(z)\leq \lambda_{D_1}(z)\leq \lambda_{D_{D_2}(z,c)}(z)=\frac{1}{\tanh c}\lambda_{D_2}(z)<K \lambda_{D_2}(z).
\]
This last inequality is (b). For (c), let $z,w\in S_{D_2}(\gamma\lvert_{[t_1,+\infty)},R)$ and let $\eta\colon[0,1]\to D_2$ be the geodesic of $D_2$ joining $z$ and $w$. By Remark~\ref{rem: convexity of sectors} the hyperbolic sector $S_{D_2}(\gamma\lvert_{[t_1,+\infty)},R)$ is a hyperbolically convex set and thus contains the geodesic $\eta$. Also, $\eta([0,1])\subset D_1$, due to part (a). So, by the domain monotonicity of the hyperbolic metric \eqref{eq:hyperbolic domain monotonicity}, we have
\[
d_{D_2}(z,w)\leq d_{D_1}(z,w)\leq\int_\eta \lambda_{D_1}(\zeta)\lvert d\zeta \rvert\leq K \int_\eta \lambda_{D_2}(\zeta)\lvert d\zeta \rvert= K\ d_{D_2}(z,w).\qedhere
\]
\end{proof}

As a corollary of Theorem \ref{thm:internal tangency hyperbolic equivalence} we show that whenever $D_1$ is internally tangent to $D_2$, the two domains share many quasi-geodesics. Before proving this result, stated in Corollary \ref{cor: quasi-geodesics in internally tangent domains} to follow, we require a corollary of the Shadowing Lemma (Theorem \ref{thm:shadowing lemma}) that can be found in \cite[Corollary 6.3.9]{BCDM-Book}.

\begin{cor}[\cite{BCDM-Book}]\label{cor:shadowing corollary}
Let $D\subsetneq\C$ be a simply connected domain and $\zeta\in\partial_CD$. If $\eta\colon[0,+\infty)\to D$ is a quasi-geodesic of $D$ landing at $\zeta$ and $\{z_n\}\subset D$, then $\{z_n\}$ converges to $\zeta$ non-tangentially in $D$ if and only if there exists a constant $R>0$, so that 
\[
d_D(z_n,\eta)\leq R, \quad \text{for all }n\in\N.
\]
\end{cor}

\begin{cor}\label{cor: quasi-geodesics in internally tangent domains}
Let $D_1\subset D_2\subsetneq \C$ be two simply connected domains such that $D_1$ is internally tangent to $D_2$ at $\zeta\in\partial_CD_2$. Also, let $\zeta_1\in\partial_CD_1$ be the prime end of $D_1$ associated to $\zeta$.
\begin{enumerate}[label={\rm (\alph*)}]
\item For any quasi-geodesic $\eta_2\colon[0,+\infty)\to D_2$ of $D_2$ landing at $\zeta$ there exists a constant $t_1\geq0$ so that  $\eta_2([t_1,+\infty))\subset D_1$, and $\eta_2\colon [t_1,+\infty)\to D_1$ is a quasi-geodesic of $D_1$. When $\eta_2([0,+\infty))\subset D_1$, $t_1$ can be chosen to be zero.
\item Any quasi-geodesic $\eta_1\colon[0,+\infty)\to D_1$ of $D_1$ landing at $\zeta_1$ is also a quasi-geodesic of $D_2$ landing at $\zeta$ in the Carath\'eodory topology of $D_2$.
\end{enumerate}
\end{cor}

\begin{proof}
For part (a), let $\eta_2\vcentcolon[0,+\infty)\to D_2$ be an $(A,B)$-quasi-geodesic of $D_2$. By Theorem \ref{thm:shadowing lemma}, there exists a geodesic $\gamma\colon[0,+\infty)\to D_2$ of $D_2$ and a constant $R>0$, so that $\eta_2(t)\in S_{D_2}(\gamma,R)$, for all $t\geq0$. For any fixed $K>1$, Theorem~\ref{thm:internal tangency hyperbolic equivalence} implies that there exists $t_0\geq 0$ such that $S_{D_2}(\gamma\lvert_{[t_0,+\infty)},R)\subset D_1$ and
\begin{equation}\label{eq:qg coro eq1}
\lambda_{D_1}(z)\leq K \lambda_{D_2}(z), \quad \text{for all}\ z\in S_{D_2}(\gamma\lvert_{[t_0,+\infty)},R).
\end{equation}
We can also find $t_1\geq0$, so that
\[
\eta_2(t)\in S_{D_2}(\gamma\lvert_{[t_0,+\infty)},R)\subset D_1,\quad \text{for all}\ t\geq t_1.
\]
This already shows that $\eta_2([t_1,+\infty))$ lies in $D_1$. Also, using \eqref{eq:qg coro eq1}, we have that for all $t_1\leq s\leq t$
\begin{align*}
\ell_{D_1}(\eta_2; [s,t])&=\int_{\eta_2\lvert_{[s,t]}}\lambda_{D_1}(z) \lvert dz \rvert\leq K \int_{\eta_2\lvert_{[s,t]}}\lambda_{D_2}(z) \lvert dz \rvert\\
&= K \ell_{D_2}(\eta_2;[s,t])\leq KA\ d_{D_2}(\eta_2(s),\eta_2(t))+KB\\
&\le KA d_{D_1}(\eta_2(s),\eta_2(t))+KB.
\end{align*}
Therefore $\eta_2\colon[t_1,+\infty)\to D_1$ is a $(KA,KB)$-quasi-geodesic of $D_1$. If, in addition, we assume that $\eta_2([0,+\infty))\subset D_1$, then setting $B'=KB+\ell_{D_1}(\eta_2;[0,t_1])$ we obtain that  $\eta_2\colon[0,+\infty)\to D_1$ is a $(KA,B')$-quasi-geodesic of $D_1$.\\
For part (b), let $\eta_1$ be a quasi-geodesic of $D_1$ landing at $\zeta_1\in\partial_CD_1$ (recall that the prime end $\zeta_1$ was constructed in Lemma \ref{lm:intersection of prime ends}). Because $D_1$ is internally tangent to $D_2$ at $\zeta\in\partial_CD_2$, there exists a geodesic $\gamma_2\colon[0,+\infty)\to D_2$ of $D_2$ landing at $\zeta$, such that $\gamma_2([0,+\infty))\subset D_1$. We claim that there exist constants $t_2\geq0$ and $R'>0$, so that $\eta_1(t)\in S_{D_2}(\gamma_2,R')$, for all $t\geq t_2$. If this were not the case, there would exist a sequence $\{t_n\}\subset[0,+\infty)$, with $t_n\xrightarrow{n\to+\infty}+\infty$, so that $\{\eta_1(t_n)\}$ is not contained in any hyperbolic sector of $D_2$ around $\gamma_2$. Observe that $\{\eta_1(t_n)\}$ converges to $\zeta_1$ non-tangentially in $D_1$ due to Corollary \ref{cor:shadowing corollary}. But, from part (a), $\gamma_2$ is also a quasi-geodesic of $D_1$. Therefore using Corollary \ref{cor:shadowing corollary}, again, along with the domain monotonicity of the hyperbolic distance \eqref{eq:hyperbolic domain monotonicity}, we can find a constant $d>0$ so that
\[
d_{D_2}(\eta_1(t_n),\gamma_2)\leq d_{D_1}(\eta_1(t_n),\gamma_2)\leq d,
\]
which implies that $\{\eta_1(t_n)\}$ is contained in the sector $S_{D_2}\left(\gamma_2,d\right)$, leading to a contradiction. Hence, we have that $\eta_1$ is eventually contained in a hyperbolic sector of $D_2$. This immediately yields that $\eta_1$ lands at $\zeta$ in the Carath\'eodory topology of $D_2$. The fact that $\eta_1$ is a quasi-geodesic of $D_2$ follows from arguments similar to those in part (a).
\end{proof}

Combining Corollary \ref{cor: quasi-geodesics in internally tangent domains} with Corollary \ref{cor:shadowing corollary} immediately yields that internally tangent domains share any sequences converging non-tangentially.

\begin{cor}\label{coro:nontangetial convergence in internally tangent domains}
    Let $D_1\subset D_2\subsetneq \C$ be two simply connected domains such that $D_1$ is internally tangent to $D_2$ at $\zeta\in\partial_CD_2$. Also, let $\zeta_1\in\partial_CD_1$ be the prime end of $D_1$ associated to $\zeta$. Any sequence $\{z_n\}\subset D_2$ converging non-tangentially to $\zeta$ in $D_2$ is eventually contained in $D_1$ and converges non-tangentially to $\zeta_1$ in $D_1$. Conversely, any sequence $\{z_n\}\subset D_1$ converging non-tangentially to $\zeta_1$ in $D_1$ also converges to $\zeta$ non-tangentially in $D_2$.
\end{cor}

Using Corollary \ref{coro:nontangetial convergence in internally tangent domains} we prove a transitivity property for internally tangent domains.

\begin{cor}\label{cor:transitivity of internal tangency}
Let $D_1\subset D_2\subset D_3\subsetneq \C$ be simply connected domains and $\zeta\in\partial_CD_3$.
\begin{enumerate}[label={\rm (\alph*)}]
\item Assume that $D_1$ is internally tangent to $D_3$ at $\zeta$. Then $D_2$ is also internally tangent to $D_3$ at $\zeta$, and $D_1$ is internally tangent to $D_2$ at the prime end of $D_2$ associated to $\zeta$.
\item If $D_2$ is internally tangent to $D_3$ at $\zeta$ and $D_1$ is internally tangent to $D_2$ at the prime end of $D_2$ associated to $\zeta$, then $D_1$ is internally tangent to $D_3$ at $\zeta$.
\end{enumerate}
\end{cor}

\begin{proof}
For part (a), the fact that $D_2$ is internally tangent to $D_3$ at $\zeta$ follows immediately from the definition. Let $\zeta_2\in\partial_CD_2$ be the prime end of $D_2$ associated to $\zeta$. We now show that $D_1$ is internally tangent to $D_2$ at $\zeta_2$. Let $\{z_n\}\subset D_2$ be a sequence converging non-tangentially to $\zeta_2$ in $D_2$. Due to Remark \ref{rem: equiv def of internal tangency}, our goal is to show that $\{z_n\}$ is eventually contained in $D_1$. Since $D_2$ is internally tangent to $D_3$ at $\zeta$, Corollary \ref{coro:nontangetial convergence in internally tangent domains} implies that $\{z_n\}$ is eventually contained in $D_3$ and converges to $\zeta$ in $D_3$. But, $D_1$ is also internally tangent to $D_3$ at $\zeta$, meaning that $\{z_n\}$ is eventually contained in $D_1$, again due to Remark \ref{rem: equiv def of internal tangency}, as required. Part (b) follows from similar arguments. 
\end{proof}

In most of our results we will consider a domain $D\subset \D$ that is internally tangent to $\D$ at a point $\zeta\in\partial\D$. In this case many of our previous statements and arguments are simpler, since the Carath\'eodory topology  of $\D$ coincides with its Euclidean topology. Furthermore, in this setting the notion of internally tangent domains is closely related to another important property of conformal maps, called ``\textit{semi-conformality}" or ``\textit{isogonality}" on the boundary. The version of this property we require is stated below and is a special case of a result by Ostrowski (see, for example, \cite[Theorem 5.5, p. 177]{GM}). We also refer to \cite{GKMR} for a recent exploration of semi-conformality. 

\begin{theorem}\label{thm:semi-conformality}
If $D\subset\D$ is a simply connected domain internally tangent to $\D$ at $\zeta\in\partial\D$, there exists a Riemann map $\phi\colon \D\to D$ so that $\angle \lim\limits_{z\to\zeta}\phi(z)=\zeta$ and 
\[
\angle \lim_{z\to\zeta}\mathrm{arg}\left(\frac{1-\overline{\zeta}\phi(z)}{1- \overline{\zeta}z}\right)\in [0,2\pi).
\] 
\end{theorem}

%%%%%%%%%%%%%%%%%%%%%%%%%%%%%%%%%%%%%%%%%%%%%%%%%%%%%%%%%
\section{Fundamental domain and semigroup-fication}\label{section:fundamental domains and semigroupfication}
%%%%%%%%%%%%%%%%%%%%%%%%%%%%%%%%%%%%%%%%%%%%%%%%%%%%%%%%%

With most of the necessary preliminary material in place, we now develop an extension of the ``semigroup-fication" technique introduced by Bracci and Roth in \cite{Bracci-Roth}. In particular, in this section we prove Theorems \ref{thm:main A} and \ref{thm:main B}.

The concept of a fundamental domain of a holomorphic map $f\colon\D\to\D$ defined in the Introduction plays a crucial role to our construction. Let us note that the condition $\D=\bigcup_{n=1}^{+\infty}f^{-n}(U)$ in the definition of a fundamental domain $U$ is equivalent to the following: for every $z_0\in \D$, there exists $n_0\in\N$ such that $f^n(z_0)\in U$, for all $n\geq n_0$.

The existence of a fundamental domain was first established by Cowen \cite[Proposition 3.1, Theorem 3.2]{Cowen}, and was based on an earlier construction by Pommerenke \cite[Theorem 2]{Pommerenke-Iteration}. He also showed that whenever $f$ is non-elliptic and $\{f^n\}$ converges non-tangentially, the fundamental domain is internally tangent to the unit disc at the Denjoy--Wolff point. A construction similar to Cowen's appears in \cite[Theorem 2.2]{CDP}. For the case of a zero-parabolic map, Contreras, D\'iaz-Madrigal and Pommerenke \cite[Theorem 5.1]{CDP2} gave an entirely different construction of a fundamental domain which is always internally tangent to the disc. All these results can be summed up in the following theorem. 

\begin{theorem}[\cites{CDP2, Cowen}]\label{thm:cowen's fundamental domain}
Let $f\colon\D\to\D$ be a non-elliptic map with Denjoy--Wolff point $\tau\in\partial\D$, and suppose that $h\colon\D\to\C$ is a Koenigs function for $f$. There exists a fundamental domain $U$ for $f$ on which $h$ is univalent. If, in addition, $f$ is hyperbolic or zero-parabolic, $U$ can be chosen to be internally tangent to $\D$ at $\tau$.
\end{theorem}

A key step in the technique developed by Bracci and Roth is a method of producing a starlike at infinity subdomain of an asymptotically starlike at infinity domain, given in \cite[Lemma 7.6]{Bracci-Roth}.

\begin{lm}[\cite{Bracci-Roth}]\label{lm:starlikefication of Koenigs domain}
Let $\Omega\subsetneq\C$ be a domain asymptotically starlike at infinity. There exists a non-empty, simply connected, starlike at infinity domain $\Omega^*\subset \Omega$ which satisfies
\[
\bigcup_{n\in\N}\left(\Omega-n\right)=\bigcup_{t\geq0}\left(\Omega^*-t\right).
\]
\end{lm}

Adopting the terminology from \cite{Bracci-Roth}, we call the subdomain $\Omega^*$ the \emph{starlike-fication} of $\Omega$. In \cite[Theorem 9.2]{Bracci-Roth} they also proved the following:

\begin{theorem}[\cite{Bracci-Roth}]\label{thm:Bracci-Roth starlikefication theorem}
    Let $\Omega\subsetneq\C$ be a simply connected domain, asymptotically starlike at infinity, such that $\bigcup_{n\in\N}\left(\Omega-n\right)$ is not a strip. Consider a Riemann map $h\colon \D\to \Omega$ satisfying $\lim_{t\to+\infty}h^{-1}(w_0+t)=\zeta\in\partial\D$, for some (any) $w_0\in\Omega$. If $\{z_n\}\subset\D$ converges non-tangentially to $\zeta$ in $\D$ then $\{h(z_n)\}$ is eventually contained in $\Omega^*$.
\end{theorem}

For the rest of this section we fix a non-elliptic map $f\colon\D\to\D$ with Denjoy--Wolff point $\tau\in\partial\D$, a Koenigs domain $\Omega$ and a Koenigs function $h\colon \D\to \Omega$. Also, suppose that $U$ is the fundamental domain for $f$ given by Theorem \ref{thm:cowen's fundamental domain}.

\begin{lm}\label{lm:h(U) asymptotically starlike}
    The domain $h(U)$ is asymptotically starlike at infinity and satisfies 
    \begin{equation}\label{eq:h(U) asymptotically starlike}
    \bigcup_{n\in\N}\left(h(U)-n\right)=\bigcup_{n\in\N}\left(\Omega-n\right).
      \end{equation}
\end{lm}

\begin{proof}
    We first show that $h(U)+1\subset h(U)$. Let $h(z)\in h(U)$, for some $z\in U$. Then, because $h$ is a Koenigs function for $f$ we have that $h(z)+1=h(f(z))$. But, $U$ is a fundamental domain for $f$, meaning that $f(z)\in U$. So, $h(f(z))\in h(U)$.\\
    Since $\Omega$ is a Koenigs domain for $f$, it is asymptotically starlike at infinity. So to complete the proof it suffices to show \eqref{eq:h(U) asymptotically starlike}. Note that we trivially have $h(U)\subset \Omega$ and so 
    \[
    \bigcup_{n\in\N}\left(h(U)-n\right)\subset\bigcup_{n\in\N}\left(\Omega-n\right).
    \]
    For the inverse inclusion, let $w\in \bigcup_{n\in\N}\left(\Omega-n\right)$. Then, $w+N\in\Omega$ for some $N\in\N$. Write $h(z)=w+N$, for some $z\in\D$. Again from the fact that $U$ is a fundamental domain for $f$, we have that $f^{n_0}(z)\in U$ for some $n_0\in\N$. So, 
    \[
    w+N+n_0=h(z)+n_0=h(f^{n_0}(z))\in h(U).
    \]
    Therefore, $w\in\left(h(U)-N-n_0\right)\subset \bigcup_{n\in\N}\left(h(U)-n\right)$.
\end{proof}

Lemma \ref{lm:h(U) asymptotically starlike} allows us to apply Lemma \ref{lm:starlikefication of Koenigs domain} on $h(U)$ in order to obtain a simply connected, starlike at infinity subdomain $h(U)^*\subset h(U)$. Moreover, $h(U)^*$ satisfies
\begin{equation}\label{eq:semigroupfication eq 1}
\bigcup_{n\in\N}\left(h(U)^*-n\right)=\bigcup_{t\geq0}\left(h(U)^*-t\right)=\bigcup_{n\in\N}\left(h(U)-n\right)=\bigcup_{n\in\N}\left(\Omega-n\right).
\end{equation}
The first equality follows from simple arguments, the second from Lemma \ref{lm:starlikefication of Koenigs domain} and the third is \eqref{eq:h(U) asymptotically starlike}.

\begin{figure}[ht]
\centering
\includegraphics[scale=0.95]{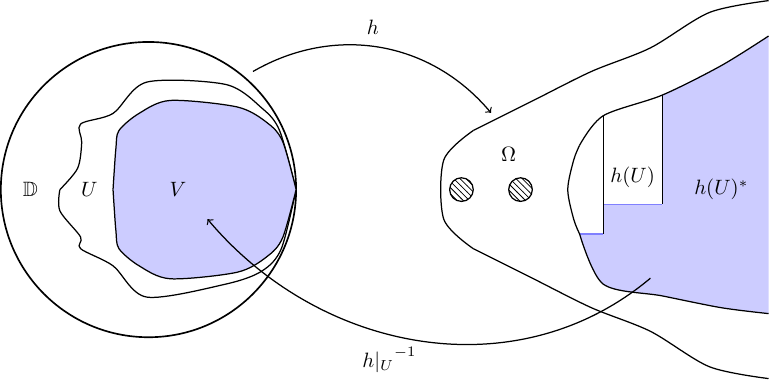}
\caption{Constructing the domain $V$.}
\label{fig: fundamental domain V} 
\end{figure}

Since $h$ is univalent on $U$, due to Theorem \ref{thm:cowen's fundamental domain}, we can define
\begin{equation}\label{eq:definition of V}
V\vcentcolon={h\lvert_{U}}^{-1}\left(h(U)^*\right),
\end{equation}
which is a simply connected subdomain of $\D$ (see Figure \ref{fig: fundamental domain V}). In fact, we are going to show that $V$ is a fundamental domain for $f$ that is internally tangent to $\D$ at $\tau$, the Denjoy--Wolff point of $f$, whenever $U$ is internally tangent to $\D$.

\begin{lm}\label{lm:V is a fundamental domain}
The set $V$ is a fundamental domain for $f$. In addition, $V$ is internally tangent to $\D$ at $\tau$ whenever $U$ has the same property.
\end{lm}

\begin{proof}
Using a translation we can assume, without loss of generality, that $0\in h(U)^*$. We first show that $V$ is a fundamental domain. Note that $V$ is simply connected by construction and $f$ is univalent on $V$ since $V\subset U$ and $U$ was a fundamental domain. Now, for any $z\in V$, we have that $h(f(z))=h(z)+1\in h(U)^*$ because $h(U)^*$ is starlike at infinity. Since $h$ is univalent on $U$, we obtain that $f(z)={h\lvert_{U}}^{-1}(h(z)+1)\in V$, showing that $f(V)\subset V$.\\
Let $z_0\in\D$. Because $U$ is a fundamental domain, there exists $n_0\in\N$, so that $f^{n_0}(z_0)\in U$. So, $h(z_0)+n_0=h(f^{n_0}(z_0))\in h(U)$, meaning that
\[
h(z_0)\in \bigcup_{n\in\N}\left(h(U)-n\right) = \bigcup_{n\in\N}\left(h(U)^*-n\right),
\]
using \eqref{eq:semigroupfication eq 1}. We conclude that there exists $n_1\in\N$ so that $h(f^{n_1}(z_0))=h(z_0)+n_1\in h(U)^*$, and the univalence of $h$ on $U$ implies that $f^{n_1}(z_0)={h\lvert_{U}}^{-1}(h(z_0)+n_1)\in V$, as required. This concludes the proof that $V$ is a fundamental domain for $f$.\\
Suppose now that $U$ is internally tangent to $\D$ at $\tau$. We split the proof in two cases depending on the type of $f$. First, we assume that $f$ is hyperbolic. Then $\bigcup_{n\in\N}\left(\Omega-n\right)$ is a horizontal strip and, up to conjugation with a translation, we can assume that
\[
\bigcup_{n\in\N}\left(h(U)^*-n\right)=\bigcup_{n\in\N}\left(\Omega-n\right)=\{x+iy\in\C\colon \lvert y \rvert < y_0\}=\vcentcolon S_0,
\]
for some $y_0\in(0,+\infty)$. Note the inclusion $h(U)^*\subset h(U)\subset S_0$. Denote by $+\infty$ the prime end at infinity of $S_0$ accessed through the positive real axis. We will show that $h(U)^*$ is internally tangent to $+\infty\in\partial_C S_0$. Because $S_0$ is symmetric with respect to the real axis, \cite[Proposition 6.1.3]{BCDM-Book} tells us that the half-line $\gamma(t)=t$, with $t\geq0$, is a hyperbolic geodesic of $S_0$ that lands at the prime end $+\infty$. Also, simple calculations show that for any $t_1\geq0$ the hyperbolic sector $S_{S_0}\left(\gamma\lvert_{[t_1,+\infty)},R\right)$ is contained in a half-strip of the form $S_{\delta,x_1}=\{x+iy\in\C\colon x_1<x, \lvert y \rvert <\delta\}$, for some $x_1\in\R$ depending on $t_1$ and $R$, and some $\delta\in(0,y_0)$ depending only on $R$. Let us fix some $R>0$ and let $x_0\in\R$ be the number for which $S_{S_0}(\gamma,R)\subset S_{\delta,x_0}$ (i.e. the $x_0$ obtained for $t=0$). Consider the vertical line segment $L=\{x+iy\in\C\colon x=x_0, \lvert y \rvert <\delta\}$. Because $L$ is compactly contained in $S_0=\bigcup_{n\in\N}\left(h(U)^*-n\right)$, there exists some $n_0$ so that $L+n_0\subset h(U)^*$. But $h(U)^*$ is starlike at infinity, so $L+n_0+t\subset h(U)^*$, for all $t>0$. This means that $S_{\delta,x_0}\cap\{z\in\C\colon \mathrm{Re}\ z \geq x_0+n_0\}\subset h(U)^*$. So, we can choose some $t_0\geq0$ large enough so that $S_{S_0}\left(\gamma\lvert_{[t_0,+\infty)},R\right)$ is contained in  $h(U)^*$, as required for $h(U)^*$ to be internally tangent to $S_0$. Now, by the transitivity property of internally tangent domains, Corollary \ref{cor:transitivity of internal tangency} (a), we get that $h(U)$ is also internally tangent to $S_0$ at $+\infty\in\partial_CS_0$.  As a slight abuse of notation, let us also denote by $+\infty\in\partial_C h(U)$ the prime end of $h(U)$ associated to $+\infty\in\partial_CS_0$. Corollary \ref{cor:transitivity of internal tangency} (a) also yields that $h(U)^*$ is internally tangent to $h(U)$ at $+\infty\in\partial_Ch(U)$. But recall that $h\lvert_U$ maps $U$ conformally onto $h(U)$ and $V$ conformally onto $h(U)^*$. So, by the conformal invariance of the notion of internally tangent domains we obtain that $V$ is internally tangent to $U$ at the prime end of $U$ associated to $\tau$. As $U$ was already internally tangent to $\D$ at $\tau$, a final use of the transitivity property of internally tangent domains, Corollary \ref{cor:transitivity of internal tangency} (b), yields that $V$ is internally tangent to $\D$ at $\tau$.\\
If $f$ is parabolic, then $\bigcup_{n\in\N}\left(\Omega-n\right)$ is not a strip and so neither is $\bigcup_{n\in\N}\left(h(U)-n\right)$, due to \eqref{eq:h(U) asymptotically starlike}. Write $\tau_U\in\partial_CU$ for the prime end of $U$ associated to $\tau$, which exists because $U$ was assumed to be internally tangent to $\D$ at $\tau$. Take a Riemann map $\phi\colon \D\to U$ with $\phi(\tau)=\tau_U$ (we identify $\phi$ with its Carath\'eodory extension, for simplicity). Let $\{z_n\}\subset \D$ be a sequence that converges non-tangentially in $\D$ to $\tau$. Due to Remark \ref{rem: equiv def of internal tangency}, in order to prove that $V$ is internally tangent to $\D$ at $\tau$, it suffices to show that $\{z_n\}$ is eventually contained in $V$. Corollary \ref{coro:nontangetial convergence in internally tangent domains} implies that $\{z_n\}$ is eventually contained in $U$ and converges non-tangentially in $U$ to $\tau_U$. Hence, by deleting finitely many terms from $\{z_n\}$ if necessary, we have that the sequence $\{\phi^{-1}(z_n)\}\subset \D$ converges to $\tau$ non-tangentially in $\D$. Now, recall that $h\lvert_U$ maps $U$ conformally onto $h(U)$, meaning that $h\lvert_U\circ \phi\colon \D\to h(U)$ is a Riemann map. It is easy to see that for any $w_0\in\Omega$, we have that $\lim_{t\to+\infty}\phi^{-1}\circ {h\lvert_U}^{-1} (w_0+t)=\tau$. Therefore, we can apply Theorem \ref{thm:Bracci-Roth starlikefication theorem} to the asymptotically starlike at infinity domain $h(U)$, the Riemann map $h\lvert_U\circ \phi$ and the sequence $\{\phi^{-1}(z_n)\}$ in order to obtain that $\{h\lvert_U\circ\phi(\phi^{-1}(z_n))\}=\{h\lvert_U(z_n)\}$ is eventually contained in $h(U)^*$. Finally, recalling that, by construction, $h\lvert_U$ maps $V$ conformally onto $h(U)^*$ yields that $\{z_n\}$ is eventually contained in $V$, as required.
\end{proof}

The material we presented so far contains the proof of Theorem \ref{thm:main A}. To be more precise, the fundamental domain $V$ we constructed satisfies all necessary properties, since $h$ is univalent on $V$, its image $h(V)=h(U)^*$ is starlike at infinity, and \eqref{eq: semigroupfication thm 1 eq1} of Theorem \ref{thm:main A} is exactly \eqref{eq:semigroupfication eq 1}. Furthermore, when $f$ is hyperbolic or zero-parabolic, $U$ is internally tangent to $\D$ at $\tau$, due to Theorem \ref{thm:cowen's fundamental domain}, and thus so is $V$, due to Lemma \ref{lm:V is a fundamental domain}. Hence, \eqref{eq: semigroupfication thm1 eq2} of Theorem \ref{thm:main A} follows immediately from Theorem \ref{thm:internal tangency hyperbolic equivalence} (a) and (b).

We now move on to Theorem \ref{thm:main B}. Define the semigroup $\phi_t\colon V\to V$ by $\phi_t(z)=h\lvert_V^{-1}(h\lvert_V(z)+t)$, for any $t\geq0$. It is easy to see that $h\lvert_V\colon V\to \C$ is a Koenigs function of $\phi_t$, with Koenigs domain $h(U)^*$. Moreover, $(\phi_t)$ is non-elliptic and \eqref{eq:semigroupfication eq 1} shows that $\phi_t$ and $f$ have the same ``type", i.e. both are either hyperbolic, zero-parabolic or positive-parabolic. By construction we have that $f\lvert_V= \phi_1$, which inductively yields that $f^n(z)=\phi_n(z)$, for any $z\in V$ and all $n\in\N$. The semigroup $(\phi_t)$ will be called the \emph{semigroup-fication of $f$ in $V$}.

These facts already prove (a) and (b) of Theorem \ref{thm:main B}.

We also note that whenever $f$ is univalent and non-elliptic, the fundamental domain $U$ in the construction we described above can be taken to be $\D$. So, $f$ can be partially embedded into a semigroup $\phi_t\colon h^{-1}\left(\Omega^*\right) \to h^{-1}\left(\Omega^*\right)$, which is conformally conjugate (by a Riemann map of $h^{-1}\left(\Omega^*\right)$) to the semigroup-fication obtained by Bracci and Roth in \cite[Section 10]{Bracci-Roth}.

Let us now discuss the convergence of the trajectories of the semigroup-fication $(\phi_t)$ of $f$. We start with two general results about semigroups. Firstly, we prove that a Lipschitz condition for the trajectories of semigroups.

\begin{lm}\label{lm:hyperbolic Lipschitz}
	Let $D\subsetneq\C$ be a simply connected domain and suppose that $(\psi_t)$ is a non-elliptic semigroup in $D$. Then, for every $z\in D$ there exists a constant $c\vcentcolon=c(z)>0$, so that 
    \[
    d_D(\psi_{t_1}(z),\psi_{t_2}(z))\le c \lvert t_1-t_2 \rvert , \quad\text{for all}\ t_1, t_2\ge 0.
    \]
\end{lm}
\begin{proof}
	Let $g$ be a Koenigs function of $(\psi_t)$ and write $\Omega:=g(D)$ for the Koenigs domain. Recall that $g$ is univalent, and $\Omega$ is simply connected and starlike at infinity. Fix $z\in D$. By the conformal invariance of the hyperbolic distance, we have 
    \begin{equation}
        d_D(\psi_{t_1}(z),\psi_{t_2}(z))=d_\Omega(g(z)+t_1,g(z)+t_2), \quad \text{for all}\ t_2\ge t_1\ge0.
    \end{equation}
    Applying Lemma \ref{lm:distance lemma} to the horizontal line segment joining $g(z)+t_1$ and $g(z)+t_2$, we get
	$$d_\Omega(g(z)+t_1,g(z)+t_2)\le\int\limits_{t_1}^{t_2}\frac{dt}{\delta_\Omega(g(z)+t)}\ .$$
	However, the Koenigs domain of a non-elliptic semigroup is starlike at infinity, meaning that $\delta_\Omega(g(z)+t)$ is an increasing function of $t\ge0$. Thus, we have that
	$$\int\limits_{t_1}^{t_2}\frac{dt}{\delta_\Omega(g(z)+t)}\le\frac{1}{\delta_\Omega(g(z))}(t_2-t_1),\quad \textup{for all }t_2\ge t_1\ge0.$$
	Combining all of the above yields $d_D(\psi_{t_1}(z),\psi_{t_2}(z))\le\frac{1}{\delta_\Omega(g(z))}(t_2-t_1)$, for all $t_2\ge t_1\ge0$. 
\end{proof}

As a corollary of Lemma \ref{lm:hyperbolic Lipschitz} we obtain an alternative proof of a recent result obtained by the second named author and Betsakos in \cite[Corollary 6.2]{BZ}.

\begin{cor}[\cite{BZ}]\label{cor:Lipschitz}
	Let $D\subsetneq\C$ be a bounded simply connected domain and suppose that $(\phi_t)$ is a non-elliptic semigroup in $D$. Then, for every $z\in D$ there exists a constant $c\vcentcolon=c(z)>0$, so that 
    \[
    \lvert \phi_{t_1}(z)-\phi_{t_2}(z)\rvert \le c \lvert t_1-t_2 \rvert , \quad\text{for all}\ t_1, t_2\ge 0.
    \]
\end{cor}
\begin{proof}
	Fix $z\in D$ and $t_2\ge t_1\ge0$. Set $\delta\vcentcolon=\textup{diam}D\in(0,+\infty)$. Using the left-hand side inequality of Lemma \ref{lm:distance lemma}, we have
	\begin{align*}
	\notag	d_D(\phi_{t_1}(z),\phi_{t_2}(z))&\ge\frac{1}{4}\log\left(1+\frac{|\phi_{t_2}(z)-\phi_{t_1}(z)|}{\min\{\delta_D(\phi_{t_1}(z)),\delta_D(\phi_{t_2}(z))\}}\right)\\
	\notag	&\ge\frac{1}{4}\log\left(1+\frac{|\phi_{t_2}(z)-\phi_{t_1}(z)|}{\delta}\right)\\
		&\ge\frac{1}{4}\dfrac{\frac{|\phi_{t_2}(z)-\phi_{t_1}(z)|}{\delta}}{1+\frac{|\phi_{t_2}(z)-\phi_{t_1}(z)|}{\delta}},
	\end{align*}
	where the last inequality follows from the fact that $\log(1+x)\ge \frac{x}{1+x}$, for $x>-1$. Rearranging, we get that
	$$d_D(\phi_{t_1}(z),\phi_{t_2}(z))\ge\frac{1}{4}\frac{|\phi_{t_2}(z)-\phi_{t_1}(z)|}{\delta+|\phi_{t_2}(z)-\phi_{t_1}(z)|}\ge\frac{|\phi_{t_2}(z)-\phi_{t_1}(z)|}{8\delta}.$$
	Finally, by Lemma \ref{lm:hyperbolic Lipschitz}, there exists $c_0(z)>0$ so that $d_D(\phi_{t_1}(z),\phi_{t_2}(z))\le c_0(z)(t_2-t_1)$, which yields the desired inequality for the constant $c(z)\vcentcolon=8\delta\ c_0(z)$.
\end{proof}

Lemma \ref{lm:hyperbolic Lipschitz} allows us to prove that the trajectories of the semigroup-fication $(\phi_t)$ of $f$ in $V$ land at the Denjoy--Wolff point of $f$, in the Euclidean topology of $\D$.

\begin{lm}\label{lm:semigroup-fication lands at tau}
    For any $z\in V$ the curve $\eta_z\colon[0,+\infty)\to\D$ with $\eta_z(t)=\phi_t(z)$ lands at $\tau$. That is,
    \[
    \lim_{t\to+\infty}\lvert \eta_z(t)-\tau\rvert=\lim_{t\to+\infty}\lvert \phi_t(z)-\tau\rvert =0.
    \]
    If, in addition, $\{f^n\}$ converges to $\tau$ non-tangentially, then $\eta_z$ lands non-tangentially in $\D$.  
\end{lm}
\begin{proof}
    Fix some $z\in V$. Note that it suffices to show that $\lim_{n\to+\infty}\lvert\phi_{t_n}(z)-\tau\rvert=0$, for any sequence $\{t_n\}\subset [0,+\infty)$ converging to $+\infty$, and that this convergence is non-tangential whenever $\{f^n\}$ converges non-tangentially. Suppose that $\{t_n\}$ is such a sequence, and observe that $f^{\lfloor t_n \rfloor}(z)=\phi_{\lfloor t_n \rfloor}(z)$, by construction of the semigroup-fication, where $\lfloor\cdot \rfloor$ is the floor function. Thus, using the domain monotonicity of the hyperbolic distance \eqref{eq:hyperbolic domain monotonicity} and Lemma \ref{lm:hyperbolic Lipschitz} yields that
    \begin{align*}
    d_\D\left(f^{\lfloor t_n \rfloor}(z), \phi_{t_n}(z)\right) &= d_\D\left(\phi_{\lfloor t_n \rfloor}(z), \phi_{t_n}(z)\right)\leq d_V\left(\phi_{\lfloor t_n \rfloor}(z), \phi_{t_n}(z)\right)\\
    &\leq c\ \lvert \lfloor t_n\rfloor -t_n\rvert <c,
    \end{align*}
    for some constant $c>0$ depending on $z$ and for all $n\in\N$. The results now follow immediately from the convergence of $\left\{f^{\lfloor t_n \rfloor}(z)\right\}$ and Lemma \ref{lm:convergence to the the boundary in D} (a).
\end{proof}

To conclude the proof of Theorem \ref{thm:main B}, note that (c) has been proved in Lemma \ref{lm:semigroup-fication lands at tau}. We now have to prove (d). That is, we have to show that $\{f^n\}$ converges to $\tau$ non-tangentially if and only if the trajectory $(\phi_t(z))$ lands at $\tau$ non-tangentially in $\D$, for any $z\in V$. The forward implication also follows from Lemma \ref{lm:semigroup-fication lands at tau}. For the converse, assume that $(\phi_t(z))$ lands at $\tau$ non-tangentially in $\D$, for any $z\in V$. Then, since $\{f^n(z)\}\subset \{\phi_t(z)\colon t\geq0\}$, for any $z\in V$ (part (a) of Theorem \ref{thm:main B}), we immediately have that $\{f^n(z)\}$ converges non-tangentially.

To conclude this section, suppose that $\{f^n\}$ converges non-tangentially. Then, $f$ is either hyperbolic or zero-parabolic, meaning that the fundamental domain $V$ is internally tangent to $\D$ at $\tau$. Write $\tau_V\in\partial_CV$ for the prime end of $V$ associated to $\tau$. According to Lemma \ref{lm:semigroup-fication lands at tau} $(\phi_t(z))$ lands at $\tau$ non-tangentially in $\D$, for any $z\in V$. Thus, using Corollary \ref{coro:nontangetial convergence in internally tangent domains} we can show that $(\phi_t(z))$ lands at $\tau_V$ non-tangentially in $V$. We just proved the following lemma. 

\begin{lm}\label{lm:DW prime end of semigroup-fication}
    If $\{f^n\}$ converges to $\tau$ non-tangentially, then the Denjoy--Wolff prime end of $(\phi_t)$ is $\tau_V\in\partial_CV$ and $(\phi_t)$ converges to $\tau_V$ non-tangentially in $V$. 
\end{lm}

%%%%%%%%%%%%%%%%%%%%%%%%%%%%%%%%%%%%%%%%%%
\section{Embedding orbits into trajectories}\label{sect:slope}
%%%%%%%%%%%%%%%%%%%%%%%%%%%%%%%%%%%%%%%%%%
With the semigroup-fication now in place, we can proceed with the proofs of Theorem \ref{thm:main slope} and its corollary, Corollary \ref{coro:main slope coro}.

We first record an immediate corollary of Theorem \ref{thm:semi-conformality}, which states that the slope of a sequence or curve remains unchanged when considered through a domain internally tangent to $\D$.

\begin{cor}\label{cor:angles in internally tangent domains}
Let $D\subset\D$ be a simply connected domain that is internally tangent to $\D$ at $\zeta\in\partial\D$. 
\begin{enumerate}[label=\textup{(\alph*)}]
\item If $\{z_n\}\subset D$ is a sequence converging to $\zeta$ with $\mathrm{Slope}_\D(z_n)\subset(-\tfrac{\pi}{2},\tfrac{\pi}{2})$, then $\mathrm{Slope}_\D(z_n) = \mathrm{Slope}_D(z_n)$. 
\item If $\gamma\colon[0,+\infty)\to D$ is a smooth curve landing at $\zeta$ and $\mathrm{Slope}_\D(\gamma)\subset(-\tfrac{\pi}{2},\tfrac{\pi}{2})$, then $\mathrm{Slope}_\D(\gamma) =\mathrm{Slope}_D(\gamma)$.
\end{enumerate}
\end{cor}

Furthermore, we need a result from the theory of semigroups. In \cite[Theorem 1.2]{BCDMGZ} the authors prove that for semigroups in $\D$, trajectories land at $\tau$ non-tangentially if and only if they are quasi-geodesics of $\D$. The conformally invariant nature of non-tangential convergence and quasi-geodesics, we can transfer this result to any simply connected domain $D\subsetneq\C$, as stated below.

\begin{theorem}[{\cite{BCDMGZ}}]\label{thm:trajectories are quasi-geodesics}
    Let $(\phi_t)$ be a non-elliptic semigroup in a simply connected domain $D\subsetneq\C$ with Denjoy--Wolff prime end $\tau\in\partial_C D$. Fix $z\in D$. Then, the trajectory $(\phi_t(z))$ lands non-tangentially at $\tau$ if and only if $(\phi_t(z))$ is a hyperbolic quasi-geodesic of $D$.
\end{theorem}

We are now ready to prove Theorem \ref{thm:main slope}. 

\begin{proof}[Proof of Theorem \ref{thm:main slope}]
	Recall that by the construction of $(\phi_t)$ in Section \ref{section:fundamental domains and semigroupfication}, since $f$ is non-elliptic, $(\phi_t)$ is also non-elliptic. Also, if $h\colon \D\to \Omega$ is the Koenigs function for $f$ used in the construction of $V$, then $h$ is univalent on $V$ and $h\lvert_V$ is a Koenigs function for the semigroup-fication $(\phi_t)$. Now, fix $z\in\D$. Since $V$ is a fundamental domain, there exists $n_0\in\N$ such that $f^n(z)\in V$, for every $n\ge n_0$. Thus, we may consider the curve $\eta_z:[0,+\infty)\to \D$ with $\eta_z(t)=\phi_t(f^{n_0}(z))$. We are going to show that $\eta_z$ has the desired properties.\\
    Firstly, from Theorem \ref{thm:main B} (a) we have that $\phi_n(z)=f^n(z)$, for all $z\in V$ and all $n\in\N$. Hence, for $n\ge n_0$,  $f^n(z)=f^{n-n_0}(f^{n_0}(z))=\phi_{n-n_0}(f^{n_0}(z))=\eta_z(n-n_0)$ and (a) is satisfied. As $\eta_z$ is a trajectory of the semigroup $(\phi_t)$, it lands at $\tau$ in the Euclidean topology of $\D$ (Lemma \ref{lm:semigroup-fication lands at tau}). Because $V\subset\D$ is bounded, Corollary \ref{cor:Lipschitz} tells us that $\eta_z$ is a Lipschitz curve. Furthermore, as $h$ is a Koenigs function, we have that 
    \begin{align}\label{eq:slope theorem eq0}
        h\lvert_V(f(\eta_z(t)))&=h\lvert_V(\eta_z(t))+1=h\lvert_V(\phi_t(f^{n_0}(z)))+1 \nonumber\\
            &=h\lvert_V(\phi_{t+1}(f^{n_0}(z)))=h\lvert_V(\eta_z(t+1)).
    \end{align}
    Using the univalence of $h\lvert_V$ in \eqref{eq:slope theorem eq0} yields 
    \begin{equation}\label{eq:slope theorem invariance}
    f(\eta_z(t))=\eta_z(t+1), \quad \textup{for all }t\ge0,
    \end{equation} 
    which immediately implies (b).\\	
	We now move on to condition (c). In order to prove that, we will first show that
    \begin{equation}\label{eq:slope theorem 2}
	\textup{Slope}_\D(\eta_z)=\cup_{w\in\eta_z}\textup{Slope}_\D(f^n(w)),
	\end{equation}
    where we use $w\in\eta_z$ to abbreviate $w\in\eta_z([0,+\infty))$. Because of (b), the inclusion $\cup_{w\in\eta_z}\textup{Slope}_\D(f^n(w))\subseteq\textup{Slope}_\D(\eta_z)$ holds trivially. For the reverse inclusion, let $s\in\textup{Slope}_\D(\eta_z)$. By definition, there exists a strictly increasing sequence $\{t_n\}\subset[0,+\infty)$ with $\lim_{n\to+\infty}t_n=+\infty$ satisfying
	\[
    \lim\limits_{n\to+\infty}\arg(1-\bar{\tau}\eta_z(t_n))=s.
    \]
	Consider the sequence $\{x_n\}\subset[0,1)$ with $x_n:=t_n-\lfloor t_n\rfloor$. Potentially taking a subsequence, we may assume that $\lim_{n\to+\infty}x_n=x_0\in[0,1]$. Write $z_0=\eta_z(x_0)\in V$. Then $f^{\lfloor t_n\rfloor}(z_0)=f^{\lfloor t_n\rfloor}(\eta_z(x_0))=\eta_z(x_0+\lfloor t_n\rfloor )$ by an inductive use of \eqref{eq:slope theorem invariance}. Applying the domain monotonicity property of the hyperbolic distance \eqref{eq:hyperbolic domain monotonicity} and Lemma \ref{lm:hyperbolic Lipschitz}, we get
	\begin{align}\label{eq:slope theorem 1}
	d_\D(f^{\lfloor t_n\rfloor}(z_0),\eta_z(t_n))&\le d_V(f^{\lfloor t_n\rfloor}(z_0),\eta_z(t_n))=d_V(\eta_z(x_0+\lfloor t_n\rfloor),\eta_z(t_n))\nonumber\\
		&=d_V\left(\phi_{x_0+\lfloor t_n\rfloor}(f^{n_0}(z)),\phi_{t_n}(f^{n_0}(z))\right)\le c|x_0+\lfloor t_n \rfloor-t_n|,
	\end{align}
	for some positive constant $c$ depending on $f^{n_0}(z)$ i.e. depending only on the point $z$ that was fixed initially. Using the convergence of $\{x_n\}$ on inequality \eqref{eq:slope theorem 1} implies that $\lim_{n\to+\infty}d_\D(f^{\lfloor t_n\rfloor}(z_0),\eta_z(t_n))=0$. So, Lemma \ref{lm:convergence to the the boundary in D} (b) is applicable and yields that $s$ is also an accumulation point of $\{\arg(1-\bar{\tau}f^{\lfloor t_n\rfloor}(z_0))\}$ which means that $s\in\cup_{w\in\eta_z}\textup{Slope}(f^n(w))$, as required.\\    
	Having established \eqref{eq:slope theorem 2}, we may proceed to the final step of the proof of (c). We distinguish three cases depending on the type of $f$.\\	
	If $f$ is positive-parabolic then either $\textup{Slope}_\D(f^n(w))=\{-\tfrac{\pi}{2}\}$ for all $w\in\D$, or $\textup{Slope}(f^n(w))=\{\tfrac{\pi}{2}\}$ for all $w\in\D$. In any case $\cup_{w\in\eta_z}\textup{Slope}_\D(f^n(w))$ is a singleton, which by \eqref{eq:slope theorem 2} leads to condition (c).\\	
	If $f$ is zero-parabolic then, as we mentioned in Section \ref{sect:dynamics}, we have that 
    \begin{equation}
        \textup{Slope}_\D(f^n(w_1))=\textup{Slope}_\D(f^n(w_2)),\quad \text{for all}\ w_1,w_2\in\D.
    \end{equation}
    Thus we may write
    \[
    \cup_{w\in\eta_z}\textup{Slope}_\D(f^n(w))=\textup{Slope}_\D(f^n(\eta_z(0)))=\textup{Slope}_\D(f^n(f^{n_0}(z)))=\textup{Slope}_\D(f^n(z)).
    \]
    Therefore, (c) is a direct consequence of \eqref{eq:slope theorem 2}.\\	
	Finally, in the case where $f$ is hyperbolic, the semigroup-fication $(\phi_t)$ is also hyperbolic. Hence, if $\eta_w(t)=\phi_t(w)$ is a trajectory of $(\phi_t)$, for some $w\in V$, then $\mathrm{Slope}_V(\eta_w)$ is a singleton contained in $(-\tfrac{\pi}{2},\tfrac{\pi}{2})$. However, by Corollary \ref{cor:angles in internally tangent domains} we know that in this case $\mathrm{Slope}_V(\eta_w)=\mathrm{Slope}_\D(\eta_w)$. Therefore, $\textup{Slope}_\D(\eta_z)$ is again a singleton, say $\{\theta\}$. By \eqref{eq:slope theorem 2} we get that $\cup_{w\in\eta_z}\textup{Slope}_\D(f^n(w))=\{\theta\}$ which in turn leads to $\textup{Slope}_\D(f^n(\eta_z(0)))=\textup{Slope}_\D(f^n(f^{n_0}(z)))=\textup{Slope}_\D(f^n(z))=\{\theta\}$ and condition (c) is proved.\\
    To conclude the proof of the theorem, suppose that $\{f^n(z)\}$ converges to $\tau$ non-tangentially in $\D$. We are going to show that the curve $\eta_z=(\phi_t(f^{n_0}(z)))$ we constructed is a quasi-geodesic of $\D$. Observe that $f$ has to be either hyperbolic or zero-parabolic. In any case the fundamental domain $V$ we constructed in Section \ref{section:fundamental domains and semigroupfication} is internally tangent to $\D$ at $\tau$. If $\tau_V\in\partial_CV$ is the prime end of $V$ associated to $\tau$, then Lemma \ref{lm:DW prime end of semigroup-fication} implies that $\tau_V$ is the Denjoy--Wolff prime end of $(\phi_t)$ and $(\phi_t)$ converges to $\tau_V$ non-tangentially in $V$. Thus any trajectory of $(\phi_t)$, and so the curve $\eta_z$, is a quasi-geodesic of $V$ landing at $\tau_V$, due to Theorem \ref{thm:trajectories are quasi-geodesics}. Using Corollary \ref{cor: quasi-geodesics in internally tangent domains} yields that $\eta_z$ is also a quasi-geodesic of $\D$. Conversely, if $\eta_z=(\phi_t(f^{n_0}(z)))$ is a quasi-geodesic of $\D$, then it necessarily lands at $\tau$ non-tangentially. By condition (c), $\{f^n(z)\}$ converges to $\tau$ non-tangentially as well.
\end{proof}

Before exploring the ramifications of Theorem \ref{thm:main slope} to the orbits of our self-map, we provide an immediate corollary of Theorem \ref{thm:main slope} concerning semigroups in $\D$. 

\begin{cor}\label{cor:slope}
	Let $(\phi_t)$ be a semigroup in $\D$. Fix $z\in\D$ and consider the trajectory $\eta_z\colon[0,+\infty)\to\D$, with $\eta_z(t)=\phi_t(z)$. Then $\textup{Slope}_\D(\eta_z)=\textup{Slope}_\D(\phi_{t_0}^n(z))$, for any $t_0>0$.
\end{cor}

This result can certainly be obtained by arguments far simpler than the ones used in Theorem \ref{thm:main slope}. We record it here, however, since to the best of our knowledge it does not appear in the literature.

Let us now consider a non-elliptic self-map of $\D$ that converges to its Denjoy--Wolff point non-tangentially. The fact that the orbits of such a function can always be embedded in quasi-geodesics of $\D$ seems to imply that they should approach the Denjoy--Wolff point in a ``controlled" manner.

To explore this idea we first prove a characterisation for sequences that behave like quasi-geodesic curves in any planar hyperbolic domain, not necessarily simply connected (or even $\delta$-Gromov hyperbolic for that matter). 

\begin{lm}\label{lm:quasi-lemma}
    Let $D$ be a hyperbolic domain and $\{z_n\}_{n=0}^{+\infty}$ a sequence in $D$, such that the sequence $\{d_D(z_n,z_{n+1})\}$ is bounded and $\lim_{n\to+\infty}d_D(z_0,z_n)=+\infty$. Then the following are equivalent:
    \begin{enumerate}[label=\rm (\alph*)]
    \item There exist constants $A\geq1$ and $B\geq0$ so that for any integers $0 \le n<m$, we have 
    \[
    \sum_{k=n}^{m-1}d_D(z_k,z_{k+1}) \leq Ad_D(z_n,z_m)+B.
    \]
    \item There exists a quasi-geodesic $\gamma\colon[0,+\infty)\to D$ and a sequence $\{t_n\}\subset[0,+\infty)$ increasing to $+\infty$, such that $z_n=\gamma(t_n)$, for all $n\in\N\cup\{0\}$.
    \end{enumerate}
\end{lm}
\begin{proof}
    Assume that condition (a) holds. We are going to construct the desired quasi-geodesic $\gamma$. For $n\in\N\cup\{0\}$, let $\gamma_n:[0,1]\to D$ be a minimal geodesic of $D$ with $\gamma_n(0)=z_n$ and $\gamma_n(1)=z_{n+1}$. That is 
    \begin{equation}\label{eq:quasi-lemma minimal}
        \ell_D(\gamma_n;[s,t])=d_D(\gamma_n(t),\gamma_n(s)),\quad \text{for any}\ 0\leq s<t\leq 1,\ \text{and all}\ n\in\N.
    \end{equation}
    Then, consider $\gamma:[0,+\infty)\to D$ to be the curve defined by $\gamma(t)=\gamma_{\lfloor t\rfloor}(t-\lfloor t\rfloor)$. It is easy to see that $\gamma(n)=\gamma_n(0)=z_n$ for all $n\in\N\cup\{0\}$, so all that remains to be proven is that $\gamma$ is a quasi-geodesic of $D$. Fix $t_1,t_2\in[0,+\infty)$ with $t_1\le t_2$. Let $n$ be the largest integer with $n\le t_1$ and $m$ the smallest integer with $m\ge t_2$. If $n+1=m$, we immediately get that 
    \begin{equation}\label{eq:quasi-geodesic sequence 0}
    \ell_D(\gamma;[t_1,t_2])=\ell_D(\gamma_n;[t_1-n,t_2-n])=d_D(\gamma_n(t_1),\gamma_n(t_2))=d_D(\gamma(t_1),\gamma(t_2)),
    \end{equation}
    due to \eqref{eq:quasi-lemma minimal} and the definition of $\gamma$. Thus, we assume that $n+1<m$. Then, because $n\leq t_1\leq t_2\leq m$, we have that
    \begin{align}\label{eq:quasi-geodesic sequence 1}
    \notag    \ell_D(\gamma;[t_1,t_2]) &\le \ell_D(\gamma;[n,m])=\sum\limits_{k=n}^{m-1}\ell_D(\gamma;[k,k+1])=\sum\limits_{k=n}^{m-1}\ell_D(\gamma_k;[0,1]) \\
    \notag    &= \sum\limits_{k=n}^{m-1}d_D(\gamma_k(0),\gamma_k(1))=\sum\limits_{k=n}^{m-1}d_D(z_k,z_{k+1})\\
    \notag    &\le A d_D(z_n,z_m) + B = A d_D(\gamma(n),\gamma(m))+B \\
        &\le A d_D(\gamma(t_1),\gamma(t_2)) +Ad_D(\gamma(n),\gamma(t_1)) +A d_D(\gamma(t_2),\gamma(m)) +B,
    \end{align}
    where we have used \eqref{eq:quasi-lemma minimal} and condition (a) for $\{z_n\}$. In addition, since the sequence $\{d_D(z_n,z_{n+1})\}$ is bounded by assumption, we have that $M=\sup_{n\ge0}d_D(z_n,z_{n+1})$ satisfies $M\in[0,+\infty)$. Now, because both points $\gamma(n)$ and $\gamma(t_1)$ belong to the minimal geodesic $\gamma_n$ we have that 
    \begin{equation}\label{eq:quasi-geodesic sequence +}
    d_D(\gamma(n),\gamma(t_1))\le d_D(\gamma(n),\gamma(n+1))=d_D(z_n,z_{n+1})<M.
    \end{equation}
    Similarly, $\gamma(t_2)$ and $\gamma(m)$ belong to the minimal geodesic $\gamma_{m-1}$, and so
    \begin{equation}\label{eq:quasi-geodesic sequence ++}
        d_D(\gamma(t_2),\gamma(m))\le d_D(\gamma(m-1),\gamma(m))=d_D(z_{m-1},z_m)<M.
    \end{equation}
    Applying \eqref{eq:quasi-geodesic sequence +} and \eqref{eq:quasi-geodesic sequence ++} to \eqref{eq:quasi-geodesic sequence 1} we obtain
    \begin{equation}\label{eq:quasi-geodesic sequence 2}
        \ell_D(\gamma;[t_1,t_2]) \le A d_D(\gamma(t_1),\gamma(t_2)) +2AM+B.
    \end{equation}
    Note that \eqref{eq:quasi-geodesic sequence 2} holds trivially even when $n+1=m$ due to \eqref{eq:quasi-geodesic sequence 0}. Thus $\gamma$ is a $(A,B')-$quasi-geodesic of $D$, where $B'=2AM+B$.\\
    For the converse, assume that $\gamma$ is a quasi-geodesic with the properties stated in (b). Then, there exist constants $A\ge1$ and $B\ge0$ such that 
    \begin{equation}\label{eq:quasi-geodesic sequence 3}
        \ell_D(\gamma;[s_1,s_2]) \le Ad_D(\gamma(s_1),\gamma(s_2)) +B,
    \end{equation}
    for all $0\le s_1\le s_2$. Fix $n,m\in\N$ with $n<m$. Then
    \begin{align*}
        \sum\limits_{k=n}^{m-1}d_D(z_k,z_{k+1})&= \sum\limits_{k=n}^{m-1}d_D(\gamma(t_k),\gamma(t_{k+1})) \le \sum\limits_{k=n}^{m-1}\ell_D(\gamma;[t_k,t_{k+1}]) \\
        &= \ell_D(\gamma;[t_n,t_m]) \le Ad_D(\gamma(t_n),\gamma(t_m))+B\\
        &= Ad_D(z_n,z_m) +B,
    \end{align*}
    which is exactly condition (a).
\end{proof}

Using Lemma \ref{lm:quasi-lemma} and Theorem \ref{thm:main slope} we can prove Corollary \ref{coro:main slope coro}.

\begin{proof}[Proof of Corollary \ref{coro:main slope coro}]
    Let $\tau\in\partial\D$ be the Denjoy--Wolff point of a non-elliptic map $f\colon\D\to\D$. Suppose that condition (a) holds and fix $z\in\D$. By the Denjoy--Wolff Theorem, we have that $\lim_{n\to+\infty}d_\D(f^0(z),f^n(z))=+\infty$. Also, by the Schwarz--Pick Lemma \eqref{eq:hyperbolic contraction}, the sequence $\{d_{\D}(f^n(z),f^{n+1}(z))\}$ is bounded above by $d_{\D}(z,f(z))$. Thus Lemma \ref{lm:quasi-lemma} is applicable to the sequence $\{f^n(z)\}$ and implies that there exists a quasi-geodesic $\gamma:[0,+\infty)\to\D$ of $\D$, such that $\{f^n(z)\}\subset \gamma([0,+\infty))$. By Theorem \ref{thm:shadowing lemma}, $\gamma$ lands at $\tau$ non-tangentially, and thus $\{f^n(z)\}$ converges to $\tau$ non-tangentially.\\
    Conversely, suppose that $\{f^n(z)\}$ converges to $\tau$ non-tangentially, for some (and hence for all) $z\in\D$. By Theorem \ref{thm:main slope} we have that there exists some $n_0\in\N$ such that the curve $\eta_z\colon[0,+\infty)\to \D$ with $\eta_z(t)=\phi_t(f^{n_0}(z))$ lands at $\tau$ and satisfies $\eta_z(0)=f^{n_0}(z)$ and $f(\eta_z(t))=\eta_z(t+1)$, for all $t\in[0,+\infty)$ (the latter condition is \eqref{eq:slope theorem invariance} in the proof of Theorem \ref{thm:main slope}). Moreover, $\eta_z$ is a quasi-geodesic of $\D$. Note that these properties of $\eta_z$ imply that $f^{n+n_0}(z)=\eta_z(n)$, for all $n\in\N$. Using the Denjoy--Wolff Theorem and the Schwarz--Pick Lemma, again, we have that Lemma \ref{lm:quasi-lemma} is applicable to the sequence $\{f^{n+n_0}(z)\}_{n=0}^{+\infty}$. So, we can find constants $A\geq1$ and $B'\geq0$ such that for all $0\leq n<m$ 
    \[
    \sum_{k=n}^{m-1}d_\D(f^{k+1+n_0}(z),f^{k+n_0}(z))\leq A\ d_\D(f^n(z),f^m(z))+B'.
    \]
    Setting 
    \[
    B=B'+\sum_{k=0}^{n_0}d_\D(f^{k+1}(z),f^k(z))+A\cdot \max_{n,m\leq n_0}\left\{d_\D(f^n(z),f^m(z))\right\},
    \]
    we obtain that for all $0\leq n<m$
    \[
    \sum_{k=n}^{m-1}d_\D(f^{k+1}(z),f^{k}(z))\leq A\ d_\D(f^n(z),f^m(z))+B,
    \]
    which is exactly \eqref{eq:quasi-lemma eq}.
\end{proof}

%%%%%%%%%%%%%%%%%%%%%%%%%%%%%%%%%%%%%%%%%%%%%%%%%%%%%%%%
\section{Rates of convergence}\label{section:convergence rates}
%%%%%%%%%%%%%%%%%%%%%%%%%%%%%%%%%%%%%%%%%%%%%%%%%%%%%%%%

In this section we examine the rate of convergence of a non-elliptic orbit. Our main goal is to prove Theorem \ref{thm:main rates} and establish several of its corollaries.

We start with a brief rundown of the results that are already known. First of all, whenever $f$ is hyperbolic, an inductive use of Julia's Lemma \eqref{eq:Julia's lemma}, along with simple arguments, yields that for every $z\in\D$ there exists a positive constant $c\vcentcolon=c(z)$ so that
\begin{equation}\label{eq:hyperbolic rate}
    \lvert f^n(z)-\tau\rvert\le c\left(f'(\tau)\right)^n,\quad \text{for all}\ n\in\N.
\end{equation}
For a proof, see \cite[Proposition 3.1]{CZZ}. Note that \eqref{eq:hyperbolic rate} implies (see \cite[Theorem 7.1]{CZZ} for details) that for every $z\in\D$ there exists a constant $c\vcentcolon=c(z)$ so that
\begin{equation}\label{eq:hyperbolic hyperbolic rate}
    d_\D(z,f^n(z))\geq \frac{n}{2}\log\frac{1}{f'(\tau)} +c, \quad \text{for all}\ n\in\N.
\end{equation}

For the case where $f$ is positive-parabolic, the authors of \cite[Theorem 7.2]{BCDM-Rates} prove that for each $z\in\D$ there exists a positive constant $c\vcentcolon=c(z)$ such that
\begin{equation}\label{eq:parabolic Euclidean rate}
 \lvert f^n(z)-\tau\rvert\le \frac{c}{n}, \quad n\in\N.
\end{equation}
In \cite[Remark 7.3]{BCDM-Rates} it is mentioned that inequality \eqref{eq:parabolic Euclidean rate} is true for univalent $f$, but a minor modification of their arguments shows that it holds in general (see \cite[Proposition 3.4]{CZZ} for a proof). Inequality \eqref{eq:parabolic Euclidean rate} yields
\[
d_\D(z,f^n(z))\geq \log n +c, \quad\text{for all } n\in\N,
\]
for some constant $c\vcentcolon =c(z)$; see \cite[Theorem 7.4]{CZZ}.

As we can see, our main contribution to the topic of the rates of convergence is for zero-parabolic self-maps of the unit disc. So, the reader can safely assume that all functions we deal with in this section are of this type.

We start our analysis with the special case where orbits converge non-tangentially. This restriction allows us to use the material of Section \ref{sect:internal tangency} in order to relate the rates of convergence of $f$ and its semigroup-fication. Note that this result contains no assumptions on the Koenigs domain of $f$. 

\begin{theorem}\label{thm:rate-zerostep-polar}
    Let $f:\D\to\D$ be a non-elliptic map with Denjoy--Wolff point $\tau\in\partial\D$. If $\{f^n\}$ converges to $\tau$ non-tangentially, then for every $z\in\D$ and every $\epsilon>0$, there exists a constant $c\vcentcolon=c(z,\epsilon)$ such that
    \begin{equation*}
    d_{\D}(z,f^n(z))\ge \dfrac{1}{4+\epsilon}\log n + c, \quad\text{for all }n\in\N.
    \end{equation*}
\end{theorem}
\begin{proof}
    Let $(\phi_t)$ be the semigroup-fication of $f$ in $V$. Since $\{f^n\}$ converges non-tangentially, $f$ is either hyperbolic or zero-parabolic and the same is true for its semigroup-fication. In either case, the fundamental domain $V$ is internally tangent to $\D$ at $\tau$. Fix $z\in\D$ and let $n_0\in\N$ be the smallest positive integer such that $f^n(z)\in V$, for every $n\ge n_0$. The fact that $\{f^n(z)\}$ converges to $\tau$ non-tangentially implies that the trajectory $\left(\phi_t\left(f^{n_0}(z)\right)\right)$ also converges to $\tau$ non-tangentially in $\D$ (see Lemma \ref{lm:semigroup-fication lands at tau}). Therefore, there exists a geodesic $\gamma:[0,+\infty)\to\D$ of $\D$ landing at $\tau$ and some $R>0$ such that $\{f^n(z)\}\cup\{\phi_t(f^{n_0}(z)):t\ge0\}\subset S_\D(\gamma,R)$. Fix $\epsilon>0$ and let $K=1+\frac{\epsilon}{4}$. Since $V$ is internally tangent to $\D$ at $\tau$, Theorem \ref{thm:internal tangency hyperbolic equivalence} implies that there exists some $t_1\ge0$ such that $S_\D(\gamma|_{[t_1,+\infty)},R)\subset V$, and
    \begin{equation}\label{eq:polar-1}
        d_\D(w_1,w_2)\le d_V(w_1,w_2)\le Kd_\D(w_1,w_2),
    \end{equation}
    for all $w_1,w_2\in S_\D(\gamma|_{[t_1,+\infty)},R)$. For the sake of simplicity, write $z_0=f^{n_0}(z)\in V$. Then, there exists $n_1\in\N$ such that $\phi_t(z_0)\in S_\D(\gamma|_{[t_1,+\infty)},R)$, for all $t\ge n_1$. Recalling that $f^n(z_0)=\phi_n(z_0)$, we also get that $f^n(z_0)\in S_\D(\gamma|_{[t_1,+\infty)},R)$, for every $n\ge n_1$. Using the triangle inequality and \eqref{eq:polar-1}, we obtain that for all $n\ge n_1$
    \begin{eqnarray}\label{eq:polar-3}
     \notag   d_{\D}(z,f^n(z_0)) &\ge& d_{\D}(f^{n_1}(z_0),f^n(z_0))-d_{\D}(0,f^{n_1}(z_0)) -d_{\D}(0,z) \\
     \notag &\ge&\frac{1}{K} d_{V}(f^{n_1}(z_0),f^n(z_0))-d_{\D}(0,f^{n_1}(z_0))-d_{\D}(0,z)\\
     &=&\frac{1}{K} d_{V}(\phi_{n_1}(z_0),\phi_n(z_0))-d_{\D}(0,f^{n_1}(z_0))-d_{\D}(0,z).
    \end{eqnarray}
    Our goal now is to estimate the quantity $d_{V}(\phi_{n_1}(z_0),\phi_n(z_0))$ in \eqref{eq:polar-3}. First, write $\tau_V\in\partial_CV$ for the prime end of $V$ associated to $\tau$, which exists because $V$ is internally tangent to $\D$ at $\tau$ (see Lemma \ref{lm:intersection of prime ends}). Recall that $\tau_V$ is the Denjoy--Wolff prime end of the semigroup-fication $(\phi_t)$ and $(\phi_t)$ converges to $\tau_V$ non-tangentially in $V$ (see Lemma \ref{lm:DW prime end of semigroup-fication}). Let $C:\D\to V$ be a Riemann map with $C(1)=\tau_V$, where, as per usual, we have identified $C$ with its Carath\'eodory extension. Then, by defining $\psi_t\vcentcolon=C^{-1}\circ\phi_t\circ C$, we get a semigroup $(\psi_t)$ of $\D$ with Denjoy--Wolff point $1$. Let $w_0\vcentcolon=C^{-1}(z_0)\in\D$. By the conformal invariance of the hyperbolic distance and the triangle inequality, we have
    \begin{align}\label{eq:polar-4}
        d_V(\phi_{n_1}(z_0),\phi_n(z_0))&=d_\D(\psi_{n_1}(w_0),\psi_n(w_0))\nonumber \\
        &\ge d_\D(0,\psi_n(w_0))-d_\D(0,\psi_{n_1}(w_0)).
    \end{align}
    But using the formula for the hyperbolic distance in $\D$, \eqref{eq:hyperbolic distance unit disk}, and the (Euclidean) triangle inequality, we have
    \begin{equation}\label{eq:polar-5}
        d_{\D}(0,\psi_n(w_0))=\frac{1}{2}\log\frac{1+|\psi_n(w_0)|}{1-|\psi_n(w_0)|}\ge\frac{1}{2}\log\frac{1}{|\psi_n(w_0)-1|} \ge \frac{1}{4}\log n +c_0,
    \end{equation}
    for some real constant $c_0$ depending on $w_0$, where the last inequality follows from rate of convergence of $(\psi_t)$ given by Theorem \ref{thm:semigroup rates}. Combining \eqref{eq:polar-3}, \eqref{eq:polar-4} and \eqref{eq:polar-5} implies that for all $n\geq n_1$, we have
    \begin{equation*}
        d_{\D}(z,f^n(z_0)) \ge \frac{1}{4K}\log n + c_1=\frac{1}{4+\epsilon}\log n +c_1,
    \end{equation*}
    where 
    \[
    c_1=\frac{c_0}{K}-\frac{d_{\D}(0,\psi_{n_1}(w_0))}{K}-d_{\D}(0,f^{n_1}(z_0))-d_{\D}(0,z).
    \]
    As a result, we have found a constant $c_1$ such that
    \begin{equation}\label{eq:polar-6}
         d_{\D}(z,f^{n+n_0}(z))=d_{\D}(z,f^n(z_0))\ge \frac{1}{4+\epsilon}\log n + c_1,
    \end{equation}
    for all $n\ge n_1$. This is the desired inequality, but only for $n\geq n_0+n_1$. For the first $n_0+n_1-1$ terms we work as follows. Let
    \[
    c_2\vcentcolon =\min\left\{d_{\D}(z,f^n(z))-\frac{1}{4+\epsilon}\log n\colon n=1,2,\dots,n_0+n_1-1\right\}.
    \]
    Then, trivially
    \begin{equation}\label{eq:polar-7}
        d_{\D}(z,f^n(z)) \ge \frac{1}{4+\epsilon}\log n+c_2,
    \end{equation}
    for all $n=1,2,\dots, n_0+n_1-1$. Tracing back the dependencies of all the constants involved in the proof, we can see that $c_1$ and $c_2$ depend only on $z$ and $\epsilon$. Thus setting $c\vcentcolon =\min\{c_1,c_2\}$ and combining \eqref{eq:polar-6} with \eqref{eq:polar-7}, we obtain the desired rate. 
\end{proof}

Our next objective is to obtain Theorem \ref{thm:main rates} using Theorem \ref{thm:rate-zerostep-polar}. We first require a well-known estimate on split planes (see, for example, \cite[Remark 6.3]{Bracci-Speeds}).

\begin{lm}\label{lm:Koebe}
    For each $z\in K\vcentcolon=\C\setminus (-\infty,-1]$, there  exists a positive constant $c\vcentcolon=c(z)$ such that
    \begin{equation}\label{eq:hyp distance Koebe}
      \frac{1}{4}\log t -c\leq  d_K(z,z+t)\le \frac{1}{4}\log t +c, \quad \text{for all }t\ge1.
    \end{equation}

\end{lm}

Next, we show that certain distances in $\Omega_\N$ can be realised as the rate of convergence of a non-elliptic self-map of the disc.

\begin{lm}\label{lm:extremal lemma}
    There exists a point $z_0\in\D$ and a non-elliptic map $f\colon \D\to\D$, such that $\{f^n\}$ converges to the Denjoy--Wolff point of $f$ non-tangentially, and 
    \begin{equation}\label{eq:extremal lemma eq}
        d_{\Omega_\N}(1,1+n)=d_\D(z_0,f^n(z_0)), \quad \text{for all}\ n\in\N.
    \end{equation}
\end{lm}

\begin{proof}
    Let $\pi\colon \D\to \Omega_\N$ be the unique universal covering with $\pi(0)=0$ and $\pi'(0)>0$. Consider the curve $\gamma\colon [0,+\infty)\to\Omega_\N$ with $\gamma(t)=t$. Lemma \ref{lm:minimal geodesic in Omega_N} shows that $\gamma$ is a geodesic of $\Omega_\N$. Since $\Omega_\N$ is symmetric with respect to the real axis, the function $g(z)= \overline{\pi(\overline{z})}$ is a universal covering of $\Omega_\N$ that satisfies $g(0)=0$ and $g'(0)>0$. So, $\overline{\pi(\overline{z})}=\pi(z)$, for all $z\in\D$, and the geodesic $\tilde{\gamma}\colon [0,+\infty)\to\D$ of $\D$ with $\tilde{\gamma}([0,+\infty))=[0,1)$ is the unique lift of $\gamma$ starting at 0, i.e. $\pi\circ \tilde{\gamma}=\gamma$ and $\tilde{\gamma}(0)=0$. Thus, there exists some point $z_0\in(0,1)$ such that $\pi(z_0)=1$. Observe that $\tilde{\gamma}(1)=z_0$. Now, consider the holomorphic function $g\colon \Omega_\N\to\Omega_\N$ with $g(z)=z+1$. Since $g(0)=1$, $\pi(0)=0$ and $\pi(z_0)=1$, there exists a unique lift $f\colon\D\to\D$ of $g$ so that $\pi\circ f=g\circ \pi$ and $f(0)=z_0$ (see, for example, \cite[Proposition 1.6.14]{Abate}). That is, we have that $\pi(f(z))=\pi(z)+1$ for all $z\in\D$. Moreover, since $g$ has no fixed points in $\Omega_\N$, $f$ is a non-elliptic self-map of $\D$. Let us write $\tau\in\partial\D$ for the Denjoy--Wolff point of $f$.\\
    We now prove that $\{f^n(z_0)\}$ converges to $\tau$ non-tangentially. Consider the function $h(z)=\overline{f(\overline{z})},\ z\in\D$, which is a holomorphic self-map of $\D$. Then, for all $z\in\D$ we have that 
    \[
    \pi(h(z))=\pi\left(\overline{f(\overline{z})}\right)=\overline{\pi(f(\overline{z}))}=\overline{\pi(\overline{z})+1}=\pi(z)+1.
    \]
    Furthermore, $h(0)=\overline{f(\overline{0})}=\overline{z_0}=z_0$ since $z_0\in(0,1)$. Thus the uniqueness of $f$ implies that $f\equiv h$, and so $f^n(0)\in(0,1)$ for all $n\in\N$. We conclude that $\{f^n(0)\}$ converges to 1 and is contained in the geodesic $\tilde{\gamma}$, as desired.\\
    Our final task is showing \eqref{eq:extremal lemma eq}. Fix $n\in\N$. Lemma \ref{lm:minimal geodesic in Omega_N} implies that $\gamma$ is a minimal geodesic of $\Omega_\N$. So,
    \begin{equation}\label{eq:extremal lemma eq1}
        d_{\Omega_\N}(1,1+n)=\ell_{\Omega_\N}(\gamma; [1,1+n]).
    \end{equation}
    But, from the fact that $\pi$ is a local isometry for the hyperbolic metric of $\Omega_\N$ (see \eqref{eq:universal cover is a local isometry}) we have that $\ell_{\Omega_\N}(\gamma; [1,1+n])=\ell_\D(\tilde{\gamma}; [1,1+n])$. Also, $\tilde{\gamma}(1)=z_0$ and, by the arguments above, $\tilde{\gamma}(1+n)=f^n(z_0)$. Since $\tilde{\gamma}$ is a geodesic of $\D$ we obtain that
    \begin{equation}\label{eq:extremal lemma eq2}
        \ell_\D(\tilde{\gamma}; [1,1+n])=d_\D(\tilde{\gamma}(1),\tilde{\gamma}(1+n))=d_\D(z_0,f^n(z_0)).
    \end{equation}
    So, \eqref{eq:extremal lemma eq1} and \eqref{eq:extremal lemma eq2} yield \eqref{eq:extremal lemma eq}.
\end{proof}

Using Lemma \ref{lm:extremal lemma} and the rate of convergence derived in Theorem \ref{thm:rate-zerostep-polar} we prove a more general estimate in $\Omega_\N$. This will certainly prove useful later in this section, but it might also be of independent interest. 

\begin{prop}\label{prop:extremal domain}
    For every $z\in\Omega_\N$ and every $\epsilon>0$, there exist two constants $c_1\vcentcolon=c_1(z,\epsilon)\in\R$ and $c_2\vcentcolon=c_2(z)>0$ such that
    \begin{equation}\label{eq:extremal hyp estimate}
        \dfrac{1}{4+\epsilon}\log n + c_1 \le d_{\Omega_\N}(z,z+n) \le \dfrac{1}{4}\log n + c_2, \quad \text{for all }n\in\N.
    \end{equation}
\end{prop}
\begin{proof}
    Fix $z\in\Omega_\N$. Note that since $\Omega_\N$ is asymptotically starlike at infinity, the quantity $d_{\Omega_\N}(z,z+n)$ is well-defined for all $n\in\N$. Considering the slit-plane $K=\C\setminus(-\infty,-1]$, we can see that $K\subset \Omega_\N$. Moreover, there exists some $n_0\vcentcolon=n_0(z)\in\N$ such that $z+n\in K$, for all $n\ge n_0$. By the triangle inequality and the domain monotonicity of the hyperbolic distance we get that for all $n\geq n_0$
    \begin{align}\label{eq:extremal 1}
      \notag  d_{\Omega_\N}(z,z+n)\le& d_{\Omega_\N}(z,z+n_0)+d_{\Omega_\N}(z+n_0,z+n)\\
    \notag  \le& d_{\Omega_\N}(z,z+n_0)+d_K(z+n_0,z+n)\\
     =& d_{\Omega_\N}(z,z+n_0)+d_K(z+n_0,z+n_0+n-n_0).
    \end{align}
    But, distances of the form $d_K(z,z+t)$ were evaluated in Lemma \ref{lm:Koebe}, where we showed that 
    \begin{equation}\label{eq:extremal 1+}
        d_K(z+n_0,z+n_0+n-n_0)\leq \dfrac{1}{4}\log(n-n_0)+c_z,
    \end{equation}
    for some constant $c_z$ depending only on $z$, and for all $n> n_0$. 
    Combining \eqref{eq:extremal 1} and \eqref{eq:extremal 1+} yields that
    \[
    d_{\Omega_\N}(z,z+n)\leq \dfrac{1}{4}\log(n-n_0)+c_z+d_{\Omega_\N}(z,z+n_0), \quad \text{for all}\ n>n_0.
    \]
    The right-hand side inequality of \eqref{eq:extremal hyp estimate} follows by observing that for all $n\in\N$
    \[
    d_{\Omega_\N}(z,z+n)\le \frac{1}{4}\log n+c_z+\max\{d_{\Omega_\N}(z,z+m)\colon m=1,2,...,n_0\}.
    \]
    We move on to the left-hand side inequality. By Lemma \ref{lm:extremal lemma} we can find a point $z_0\in\D$ and a non-elliptic map $f\colon\D\to\D$ so that $\{f^n\}$ converges to its Denjoy--Wolff point non-tangentially, and 
    \begin{equation}\label{eq:extremal 1++}
        d_{\Omega_\N}(1,1+n)=d_\D(z_0,f^n(z_0)), \quad \text{for all } n\in\N.
    \end{equation}
    The non-tangential convergence of $\{f^n\}$ implies that Theorem \ref{thm:rate-zerostep-polar} is applicable, and so for any $\epsilon>0$ we may find a constant $c:=c(z_0,\epsilon)$ such that
    \begin{equation}\label{eq:extremal 2}
        d_{\D}(z_0,f^n(z_0))\ge \frac{1}{4+\epsilon}\log n+c, \quad\textup{for all }n\in\N.
    \end{equation}
    Combining \eqref{eq:extremal 1++} and \eqref{eq:extremal 2} yields that
    \begin{equation}\label{eq:extremal 3}
    d_{\Omega_\N}(1,1+n)\geq \frac{1}{4+\epsilon}\log n+c.
    \end{equation}
    A simple use of the triangle inequality yields 
    \[
    d_{\Omega_\N}(z,z+n)\geq d_{\Omega_\N}(1,1+n)-d_{\Omega_\N}(1,z)-d_{\Omega_\N}(1+n,z+n).
    \]
    But $d_{\Omega_\N}(1+n,z+n)\leq d_{\Omega_\N}(1,z)$ by the Schwarz--Pick Lemma applied to the holomorphic self-map of $\Omega_\N$ with $z\mapsto z+n$. Therefore, we can conclude that 
    \begin{equation}
        d_{\Omega_\N}(z,z+n)\geq d_{\Omega_\N}(1,1+n)-2d_{\Omega_\N}(1,z)\geq \frac{1}{4+\epsilon}\log n+c-2d_{\Omega_\N}(1,z),
    \end{equation}
    which is exactly the desired inequality. 
\end{proof}

Note that as an immediate consequence of Proposition \ref{prop:extremal domain} we obtain that
\begin{equation}\label{eq:discrete Koebe rate}
\lim_{n\to+\infty}\frac{d_{\Omega_\N}(z,z+n)}{\log n}=\frac{1}{4}, \quad \text{for all}\ z\in\Omega_\N.
\end{equation}
Proposition \ref{prop:extremal domain} allows us to obtain estimates on hyperbolic distances in a class of domains larger than the class of asymptotically starlike at infinity domains, as stated below. This was stated in Proposition \ref{prop:estimate on Koenigs domains}, which we now prove. 

\begin{proof}[Proof of Proposition \ref{prop:estimate on Koenigs domains}]
    Since $\Omega$ is not the whole complex plane, the property $\Omega+1\subset\Omega$ implies that $\Omega$ is a hyperbolic domain. Thus the quantity $d_\Omega(z,z+n)$ is well-defined. Moreover, there exists a translation $T(z)=z+c$, for some $c\in\C$, so that $T(\Omega)\subset\Omega_\N$. Therefore, using the conformal invariance and the domain monotonicity of the hyperbolic distance, we obtain that $d_\Omega(z,z+n)\geq d_{\Omega_\N}(z+c,z+c+n)$, for all $z\in\Omega$ and all $n\in\N$. The result now follows immediately from Proposition \ref{prop:extremal domain}.
\end{proof}

We now use Proposition \ref{prop:extremal domain} and Theorem \ref{thm:rate-zerostep-polar} to prove Theorem \ref{thm:main rates}.

\begin{proof}[Proof of Theorem \ref{thm:main rates}]
    Let $h$ be the Koenigs function of $f$ and $\Omega\subsetneq\C$ be its Koenigs domain. Using a translation, we may assume without loss of generality that $\Omega\subset\Omega_\N$. Since $h$ is a Koenigs function, we have that $h(f^n(z))=h(z)+n$, for all $z\in\D$ and all $n\in\N$. Therefore, by the Schwarz--Pick Lemma \eqref{eq:hyperbolic contraction} and the domain monotonicity \eqref{eq:hyperbolic domain monotonicity}, we deduce that
    \begin{equation*}
        d_{\D}(z,f^n(z))\ge d_{\Omega}(h(z),h(z)+n) \ge d_{\Omega_\N}(h(z),h(z)+n),
    \end{equation*}
    for all $z\in\D$ and all $n\in\N$. Proposition \ref{prop:extremal domain} now yields \eqref{eq:hyperbolic rate eq}.
\end{proof}

\begin{rem}\label{rem:sharpness of rates}
    It is known that the inequality \eqref{eq:hyperbolic rate eq} in Theorem \ref{thm:main rates} is sharp, in the sense that we can find a non-elliptic map $f\colon\D\to\D$, with 
    \begin{equation}\label{eq:sharpness of hyperbolic rate}
        \lim\limits_{n\to+\infty}\frac{d_{\D}(z,f^n(z))}{\log n}=\frac{1}{4}, \quad\text{for all }z\in\D.
    \end{equation}
    This is due to the fact that if $K=\C\setminus(-\infty,-1]$ and $g\colon \D\to K$ a Riemann map, then Lemma \ref{lm:Koebe} shows that the non-elliptic, univalent map $f\colon\D\to\D$ defined by $f(z)=g^{-1}(g(z)+1)$ satisfies \eqref{eq:sharpness of hyperbolic rate}. We note, however, that the non-elliptic map constructed in the proof of Lemma \ref{lm:extremal lemma} provides an example of a non-univalent function satisfying \eqref{eq:sharpness of hyperbolic rate}.
\end{rem}

We can now proceed to evaluating the Euclidean rate at which the iterates approach the Denjoy--Wolff point. First, note that Theorem \ref{thm:main rates 2} follows immediately from Theorem \ref{thm:main rates}, the triangle inequality and the following estimate derived by formula \eqref{eq:hyperbolic distance unit disk}
\begin{equation}\label{eq:euclidean rate}
  e^{-2d_\D(0,z)}\leq 1- \lvert z \rvert\leq  2e^{-2d_\D(0,z)}, \quad \text{for all} \ z\in\D.
\end{equation}

Next, we provide estimates on the Euclidean distance between the orbit and the Denjoy--Wolff point.

\begin{theorem}\label{thm:Euclidean rate general}
    Let $f\colon\D\to\D$ be a non-elliptic map with Denjoy--Wolff point $\tau\in\partial\D$ and whose Koenigs domain is not the whole complex plane. Then, for every $z\in\D$ and every $\epsilon>0$ there exists a positive constant $c_1\vcentcolon=c_1(z,\epsilon)$ such that
    \begin{equation}\label{eq:Euclidean rate 1}
        |f^n(z)-\tau|\le c_1\ n^{-\frac{1}{4+\epsilon}}, \quad\text{for all }n\in\N.
    \end{equation}
    If, in addition, $\{f^n\}$ converges to $\tau$ non-tangentially, then for every $z\in\D$ and every $\epsilon>0$ there exists a positive constant $c_2\vcentcolon=c_2(z,\epsilon)$ such that
    \begin{equation}\label{eq:Euclidean rate 2}
        |f^n(z)-\tau|\le c_2\ n^{-\frac{1}{2+\epsilon}}, \quad\text{for all }n\in\N.
    \end{equation}
\end{theorem}
\begin{proof}
    We start with the proof of \eqref{eq:Euclidean rate 1}. Fix $z\in\D$ and $\epsilon>0$. Since $f$ is non-elliptic and its Denjoy--Wolff point is $\tau$, Julia's Lemma yields the existence of a constant $R>0$ such that $\lvert f^n(z)-\tau\rvert^2 \le R(1-\lvert f^n(z) \rvert ^2)$, for all $n\in\N$. Note that $R$ only depends on the choice of $z$. In fact, the least possible $R$ for which this inequality holds is exactly $R=\lvert z-\tau\rvert ^2/(1-\lvert z\rvert ^2)$. Using elementary calculations, the formula for $d_\D$ in \eqref{eq:hyperbolic distance unit disk} and the triangle inequality, we find that
    \begin{align}\label{eq:rate of convergence 1}
        \lvert f^n(z)-\tau\rvert^2 &\le 4R\frac{1-\lvert f^n(z)\rvert}{1+\lvert f^n(z)\rvert} = 4Re^{-2d_{\D}(0,f^n(z))}\nonumber\\
        &\le 4R e^{2d_{\D}(0,z)}e^{-2d_{\D}(z,f^n(z))}.
    \end{align}
    Recall that by Theorem \ref{thm:main rates}, for the chosen $z$ and $\epsilon$, there exists a constant $c_0\vcentcolon=c_0(z,\epsilon)$ so that
    \begin{equation}\label{eq:rate of convergence 2}
        d_{\D}(z,f^n(z))\ge \frac{1}{4+\epsilon}\log n +c_0, \quad\textup{for all }n\in\N.
    \end{equation}
    Applying \eqref{eq:rate of convergence 2} to \eqref{eq:rate of convergence 1} implies that
    $$\lvert f^n(z)-\tau \rvert \le 2\sqrt{R}e^{d_{\D}(0,z)}e^{-c_0}\frac{1}{n^{\frac{1}{4+\epsilon}}},$$
    as desired.\\
    For the proof of \eqref{eq:Euclidean rate 2}, assume that $\{f^n\}$ converges to $\tau$ non-tangentially. As a result, for a fixed $z\in\D$, there exists a Stolz angle of $\D$ with vertex at $\tau$ which contains the orbit $\{f^n(z)\}$. To be more precise, there exists some $K>1$ so that $\lvert f^n(z)-\tau \rvert \le K(1-\lvert f^n(z)\rvert)$, for all $n\in\N$ (see \eqref{eq:Stolz}). Proceeding just like we did above, we obtain \eqref{eq:Euclidean rate 2}.
\end{proof}

The next corollary summarises most of our results in this section so far, and includes Corollary \ref{cor:main rates cor} from the Introduction. 

\begin{cor}\label{cor:rates liminf}
    Let $f:\D\to\D$ be a non-elliptic map with Denjoy--Wolff point $\tau\in\partial\D$, and whose Koenigs domain is not the whole complex plane. Then:
    \begin{enumerate}
        \item[\textup{(a)}] $\liminf\limits_{n}\frac{d_{\D}(z,f^n(z))}{\log n}\ge\frac{1}{4}$, for all $z\in\D$;
        \item[\textup{(b)}] $\limsup\limits_{n}\frac{\log|f^n(z)-\tau|}{\log n}\le -\frac{1}{4}$, for all $z\in\D$;
        \item[\textup{(c)}] if $\{f^n\}$ converges to $\tau$ non-tangentially, then $\limsup\limits_{n}\frac{\log|f^n(z)-\tau|}{\log n}\le-\frac{1}{2}$, for all $z\in\D$.
    \end{enumerate}
\end{cor}
\begin{proof}
    Let $\epsilon>0$ and fix $z\in\D$. Then, by Theorem \ref{thm:main rates}, there exists a real constant $c:=c(z,\epsilon)$ such that $d_{\D}(z,f^n(z))\ge \frac{1}{4+\epsilon}\log n + c$, for every $n\in\N$. Dividing by $\log n$ and taking limits as $n\to+\infty$, we find
    \begin{equation*}
         \liminf\limits_{n}\frac{d_{\D}(z,f^n(z))}{\log n}\ge\frac{1}{4+\epsilon}.
    \end{equation*}
   The choice of $\epsilon>0$ was arbitrary, so letting $\epsilon\to0$, we deduce (a). The proofs of statements (b) and (c) are similar, albeit with the use of Theorem \ref{thm:Euclidean rate general}.
\end{proof}

We now establish a sharper rate for the special case where the Koenigs domain $\Omega$ has non-polar boundary. Observe that this condition implies that $\Omega\subsetneq \C$, but is certainly much stronger.

\begin{theorem}\label{thm:rate-zerostep-nonpolar}
    Let $f:\D\to\D$ be a non-elliptic map with Denjoy--Wolff point $\tau$ and Koenigs domain $\Omega$. Suppose that $\partial\Omega$ is non-polar. Then, for each $z\in\D$ there exists a positive constant $c\vcentcolon=c(z)$ such that
    $$|f^n(z)-\tau|\le\frac{c}{\sqrt{n}}, \quad \text{for all } n\in\N.$$
\end{theorem}
\begin{proof}
    Let $h\colon \D\to\Omega$ be a Koenigs function for $f$ and $(\phi_t)$ the semigroup-fication of $f$ on $V$. Fix $z\in\D$. Since $V$ is a fundamental domain for $f$, there exists an $n_0\in\N$ such that $f^n(z)\in V$, for all $n\ge n_0$. Write $z_0=f^{n_0}(z)$, so that $f^n(z_0)=\phi_n(z_0)$, for every $n\in\N$. Also define the sequence of curves $\gamma_n=\{\phi_t(z_0):t\ge n\}\subset\D$ for $n\in\N$. For the rest of the proof we consider $n\in\N$ to be fixed. Observe that
    \begin{equation}\label{eq:zero-rate-1}
        \lvert f^n(z_0)-\tau\rvert=\lvert\phi_n(z_0)-\tau\rvert\le \textup{diam}[\gamma_n].
    \end{equation}
    Recall that for the semigroup-fication we have $\lim_{t\to+\infty}\lvert \phi_t(z_0)-\tau\rvert=0$ (see Lemma \ref{lm:semigroup-fication lands at tau}). So, without loss of generality we may assume that $\phi_t(z_0)\ne0$, for all $t\ge n$. As a result, we can apply Theorem \ref{thm:solynin} and formula \eqref{eq:harmonic measure arc}, to find
    \begin{equation}\label{eq:zero-rate-2}
        \omega(0,\gamma_n,\D\setminus \gamma_n)\ge \frac{1}{\pi}\arcsin\left(\frac{\textup{diam}[\gamma_n]}{2}\right).
    \end{equation}
    Combining inequalities \eqref{eq:zero-rate-1} and \eqref{eq:zero-rate-2}, we obtain $\vert f^n(z_0)-\tau\rvert\le 2\pi\omega(0,\gamma_n,\D\setminus\gamma_n)$. Since the Koenigs function $h:\D\to\Omega$ is holomorphic, the subordination principle in \eqref{eq:harmonic measure subordination} gives $\omega(0,\gamma_n,\D\setminus\gamma_n)\le\omega(h(0),h(\gamma_n),\Omega\setminus h(\gamma_n))$. Note that the harmonic measure is well-defined and non-trivial on $\Omega\setminus h(\gamma_n)$, since $\partial\Omega$ is non-polar. Using a translation, we can always assume that $h(0)=0$, for the sake of simplicity. As $(\phi_t)$ is a semigroup on $V$, each curve $\gamma_n$ is a subset of $V$. In addition, the function $h$ is univalent on $V$, and its restriction $h\lvert_V$ is a Koenigs function for $(\phi_t)$. Consequently, $h(\gamma_n)=\{h(z_0)+t:t\ge n\}=\vcentcolon\Gamma_n$. Summing up, we have
    \begin{equation}\label{eq:zero-rate-3}
        \lvert f^n(z_0)-\tau\rvert\le 2\pi\ \omega(0,\Gamma_n,\Omega\setminus \Gamma_n).
    \end{equation}
     So our goal is to estimate the harmonic measure in the right-hand side of \eqref{eq:zero-rate-3}. As $\partial\Omega$ is non-polar, we can find some $w\in\C$ and $\delta\in(0,1)$ so that the intersection $D(w,\delta)\cap(\C\setminus\Omega)$ is non-empty and non-polar. The fact that $\Omega$ is asymptotically starlike at infinity means that for each $k\in\N$, the intersection $C_k:=D(w-k,\delta)\cap(\C\setminus\Omega)$ remains non-empty and non-polar. For each $m\in\N$ consider the domain $D_m=\C\setminus\cup_{k\ge m}\overline{C_k}$. Notice that $\Omega\subset D_m$ while $D_m\subset D_{m+1}$, for all $m\in\N$. Furthermore, consider the complements of horizontal half-strips 
    $$S_m:=\C\setminus\{z\in\C:\textup{Im}z\in[\textup{Im}w-\delta,\textup{Im}w+\delta], \textup{Re}z\in(-\infty,\textup{Re}w-m+\delta]\}, \quad m\in\N.$$
    By construction, we have that $S_m\subset S_{m+1}$ and $S_m\subset D_m$, for every $m\in\N$. Finally, consider the vertical half-plane
    $$H_m:=\{z\in\C:\textup{Re}z>\textup{Re}w-m+\delta\}, \quad m\in\N.$$
    Then $\Gamma_n\subset H_m\subset S_m\subset D_m$, for all $m\in\N$ large enough. Let us note that the inclusions $H_m\subset S_m\subset D_m$ hold for all $m\in\N$, but $\Gamma_n\subset H_m$ might not hold for the first few $m$. Thus increasing $m$ and relabelling as necessary, we can assume that
    \begin{equation}\label{eq:zero-rate eq+}
        \Gamma_n\subset H_m\subset S_m\subset D_m\quad \text{and}\quad 0\in H_m, \quad \text{for all}\ m\in\N.
    \end{equation}
    For $m\in\N$ let us write $x_m:=\textup{Re}w-m+\delta$, which due to \eqref{eq:zero-rate eq+} satisfies $x_m<0$ (see Figure \ref{fig:non-polar rate} for this construction).
    \begin{figure}[ht]
    \centering
    \includegraphics[scale=1.1]{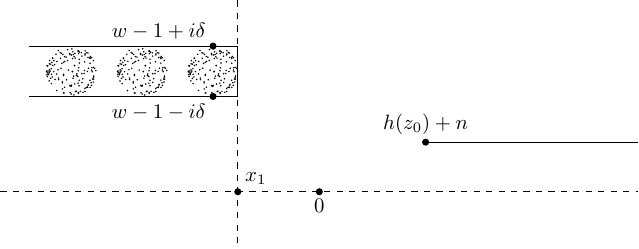}
    \caption{The construction in the proof of Theorem \ref{thm:rate-zerostep-nonpolar}}\label{fig:non-polar rate}
\end{figure}

    \noindent Now consider the M\"{o}bius transformations
    $$\Psi_m(z)=\frac{z-1-x_m}{z+1-x_m} \quad \textup{and}\quad \Phi_m(z)=\frac{z-\frac{x_m+1}{x_m-1}}{1-\frac{x_m+1}{x_m-1}z}.$$
    It is easy to check that $\Psi_m$ maps $H_m$ conformally onto the unit disc. Observe that the half-line $\Gamma_n$ is a hyperbolic geodesic of the half-plane $H_m$, emanating from $h(z_0)+n$ and landing at $\infty$. Hence, by the conformal invariance of the hyperbolic distance, $\Psi_m(\Gamma_n)$ is a geodesic in $\D$, emanating from $\Psi_m(h(z_0)+n)$ and landing at $\Psi_m(\infty)=1$. In addition, $\Phi_m$ is a conformal automorphism of the unit disk fixing $1$. As a consequence, the curve $\Phi_m\circ\Psi_m(\Gamma_n)$ is a geodesic of $\D$ emanating from $\Phi_m\circ\Psi_m(h(z_0)+n)$ and landing at $1$. Direct calculations show that $\Phi_m\circ\Psi_m(h(z_0)+n)=\frac{h(z_0)+n}{h(z_0)+n-2x_m}$, which in turn leads to 
    \begin{equation}\label{eq:zero-rate eq++}
    \lim_{m\to+\infty}\Phi_m\circ\Psi_m(h(z_0)+n)=0.
    \end{equation}
    Intuitively, \eqref{eq:zero-rate eq++} tells us that as $m$ increases to $+\infty$, the curve $\Phi_m\circ\Psi_m(\Gamma_n)$ ``transforms" into the radius of $\D$ landing at 1. We deduce that
    $$\lim\limits_{m\to+\infty}\omega(0,\Phi_m\circ\Psi_m(\Gamma_n),\D\setminus(\Phi_m\circ\Psi_m(\Gamma_n)))=1.$$
    Since M\"obius maps are homeomorphisms of the Riemann sphere, the subordination principle of the harmonic measure, \eqref{eq:harmonic measure subordination}, holds with equality. So, we get that
    $$\lim\limits_{m\to+\infty}\omega(0,\Gamma_n,H_m\setminus\Gamma_n)=1.$$
    Recall that $H_m\subset S_m\subset D_m$, for all $m\in\N$, due to \eqref{eq:zero-rate eq+}. So, the domain monotonicity property of harmonic measure, \eqref{eq:harmonic measure monotonicity}, and the fact that the harmonic measure is always bounded above by $1$, imply that
    \begin{equation}\label{eq:zero-rate-4}
        \lim\limits_{m\to+\infty}\omega(0,\Gamma_n,S_m\setminus\Gamma_n)=\lim\limits_{m\to+\infty}\omega(0,\Gamma_n,D_m\setminus\Gamma_n)=1.
    \end{equation}
    Because of \eqref{eq:zero-rate-4}, there exists $m_0\in\N$ so that
    \begin{equation}\label{eq:zero-rate-5}
        \omega(0,\Gamma_n,D_m\setminus\Gamma_n)\le 2\omega(0,\Gamma_n,S_m\setminus\Gamma_n), \quad \textup{for all }m\ge m_0.
    \end{equation}
    By the construction of the domains $D_m$, we have that $\Omega\subset D_{m_0}$. Hence, combining \eqref{eq:zero-rate-3}, \eqref{eq:zero-rate-5} with another use of the domain monotonicity of the harmonic measure, yields
    \begin{equation}\label{eq:zero-rate-6}
        \lvert f^n(z_0)-\tau\rvert\le 4\pi \ \omega(0,\Gamma_n,S_{m_0}\setminus\Gamma_n).
    \end{equation}
    For the final steps of the proof, consider the slit plane 
    \[
    D:=\C\setminus\{w-m_0+\delta-t:t\ge0\},
    \]
    and observe that $S_m\subset D$, for every $m=1,2,...,m_0$. With a final use of the monotonicity property of the harmonic measure on \eqref{eq:zero-rate-6}, we get that
    \begin{equation}\label{eq:zero-rate-7}
        \lvert f^n(z_0)-\tau\rvert\le4\pi\ \omega(0,\Gamma_n,D\setminus\Gamma_n).
    \end{equation}
    However, by \cite[Proposition 3.5]{BCDM-Rates}, there exists a positive constant $c_0$ depending on $z_0$ (and by extension on $z$) such that
    \begin{equation}\label{eq:zero-rate-8}
        \omega(0,\Gamma_n,D\setminus\Gamma_n)\le \frac{c_0}{\sqrt{n}}.
    \end{equation}
    A combination of \eqref{eq:zero-rate-7} and \eqref{eq:zero-rate-8} leads to
    \begin{equation}\label{eq:zero-rate eq+++}
         \lvert f^{n+n_0}(z)-\tau\rvert=\lvert f^n(z_0)-\tau\rvert\le \frac{4\pi c_0}{\sqrt{n}}.
    \end{equation}
   Since $n$ was chosen arbitrarily, \eqref{eq:zero-rate eq+++} is true for every $n\in\N$. Now note that since $n_0$ is fixed, we can find a constant $c_1>0$ depending on $z$ so that $\lvert f^n(z)-\tau\rvert\le c_1/\sqrt{n}$, for every $n=1,2,\dots,n_0$. Taking $c=\max\{4\pi c_0,c_1\}$ yields the desired
    \begin{equation*}
    \lvert f^n(z)-\tau\rvert\le \frac{c}{\sqrt{n}}, \quad\textup{for all }n\in\N.\qedhere
    \end{equation*}
\end{proof}

We conclude this section by examining the lower bounds for the rates of convergence to the Denjoy--Wolff point. The proof requires the following result of Arosio and Bracci \cite[Definition 2.5, Proposition 5.8]{Arosio-Bracci}.
\begin{prop}[\cite{Arosio-Bracci}]\label{prop:arosio bracci}
    Let $f\vcentcolon\D\to\D$ be a non-elliptic map with Denjoy--Wolff point $\tau$. Then
    \begin{equation*}
        \lim\limits_{n\to+\infty}\frac{d_{\D}(z,f^n(z))}{n}=-\frac{\log f'(\tau)}{2}, \quad\text{for all }z\in\D.
    \end{equation*}
\end{prop}
The original result of Arosio and Bracci is in fact valid for any non-elliptic self-map of the unit ball in $\C^n$, where $d_\D$ is replaced by the Kobayashi metric. Moreover, in \cite[Proposition 5.8]{Arosio-Bracci} the term $\frac{\log f'(\tau)}{2}$ is simply $\log f'(\tau)$. This is due to a small discrepancy in the definition of the hyperbolic metric of $\D$.

\begin{cor}\label{cor:lower bound}
    Let $f\vcentcolon\D\to\D$ be a non-elliptic map with Denjoy--Wolff point $\tau$. Then, for every $z\in\D$ and every $\epsilon\in(0,1)$ there exists a positive constant $c_0\vcentcolon=c_0(z,\epsilon)$ such that
    \begin{equation*}
        \lvert f^n(z)-\tau \rvert \ge c_0 \left(\epsilon f'(\tau)\right)^n, \quad \text{for all }n\in\N.
    \end{equation*}
\end{cor}

\begin{proof}
   Fix $z\in\D$. By some quick computations and \eqref{eq:hyperbolic distance unit disk}, we obtain
   \begin{equation}\label{eq:lower bound 1}
       \lvert f^n(z)-\tau \rvert \ge 1-\lvert f^n(z)\rvert \ge e^{-2d_{\D}(0,f^n(z))} \ge e^{-2d_{\D}(0,z)}e^{-2d_{\D}(z,f^n(z))},
   \end{equation}
   for all $n\in\N$. Fix $\epsilon\in(0,1)$. Then $-(\log\epsilon)/2>0$ and due to Proposition \ref{prop:arosio bracci}, there exists some $n_0\in\N$ such that
   \begin{equation}\label{eq:lower bound 2}
       d_{\D}(z,f^n(z))\le \left(-\frac{\log f'(\tau)}{2}-\frac{\log\epsilon}{2}\right)n, \quad\textup{for all }n\ge n_0.
   \end{equation}
   Combining \eqref{eq:lower bound 1} and \eqref{eq:lower bound 2}, we deduce
   \begin{equation*}
       \lvert f^n(z)-\tau\rvert  \ge e^{-2d_{\D}(0,z)}\left(\epsilon f'(\tau)\right)^n,
   \end{equation*}
   for all $n\ge n_0$. The result for the first $n_0-1$ terms follows by simple modifications of the constant involved. Note that this new constant will depend on the number $n_0$ which is exclusively dependent on the choice of $z$ and $\epsilon$.
\end{proof}

An analogue of Corollary \ref{cor:lower bound} is satisfied by semigroups in $\D$; cf. \cite[Theorem 2.4]{BCZZ}. In particular, the result in \cite{BCZZ} is sharp, which means that the lower bound in Corollary \ref{cor:lower bound} is also sharp.

\section{Composition operators}\label{sect:composition operator}
Here we apply our work in Section \ref{section:convergence rates} to obtain estimates for the norms of composition operators. We offer a brief rundown of the necessary background, taken from \cites{CMC, Zhu}. 

The \textit{Hardy space} $H^p$ of the unit disc, for $p\ge1$, consists of all holomorphic functions $g\vcentcolon\D\to\C$ such that
\begin{equation*}
    \sup\limits_{r\in[0,1)}\int\limits_{0}^{2\pi}\lvert g(re^{i\theta})\rvert^p d\theta <+\infty.
\end{equation*}
Each holomorphic $f\colon\D\to\D$ induces a bounded composition operator $C_f\colon H^p\to H^p$ with $C_f(g)=g\circ f$, whose norm satisfies the following estimate:

\begin{lm}[{\cite[Corollary 3.7]{CMC}}]\label{lm:growth Hardy}
    Let $f\vcentcolon\D\to\D$ be a holomorphic function and consider the composition operator $C_f\vcentcolon H^p\to H^p$, $p\ge1$. Then
    \begin{equation*}
        \left(\frac{1}{1-\lvert f(0)\rvert^2}\right)^{\frac{1}{p}} \le \rvert\rvert C_f \lvert\lvert_{H^p} \le \left(\frac{1+\lvert f(0)\rvert}{1-\lvert f(0)\rvert}\right)^{\frac{1}{p}},
    \end{equation*}
    where $\lVert\cdot \rVert_{H^p}$ denotes the operator norm in $H^p$.
\end{lm}

For $p\ge1$ and $\alpha>-1$, we consider the \textit{weighted Bergman space} $A^p_\alpha$ of the unit disc, which consists of all holomorphic functions $g\vcentcolon\D\to\C$ such that
\begin{equation*}
    \int\limits_{\D}\lvert g(z) \rvert ^p (1-\lvert z\rvert^2)^\alpha dA(z)<+\infty,
\end{equation*}
where $dA$ denotes the normalized Lebesgue area measure. Composition operators $C_f\colon A^p_\alpha\to A^p_\alpha$ are defined similarly to the case of $H^p$, are bounded and their norm satisfies:

\begin{lm}[{\cite[Section 11.3]{Zhu}}]\label{lm:growth Bergman}
    Let $f\vcentcolon\D\to\D$ be a holomorphic function and consider the composition operator $C_f\vcentcolon A^p_\alpha\to A^p_\alpha$, $p\ge1$, $\alpha>-1$. Then
    \begin{equation*}
        \left(\frac{1}{1-\lvert f(0)\rvert^2}\right)^{\frac{2+\alpha}{p}} \le \rvert\rvert C_f \lvert\lvert_{A^p_\alpha} \le \left(\frac{1+\lvert f(0)\rvert}{1-\lvert f(0)\rvert}\right)^{\frac{2+\alpha}{p}},
    \end{equation*}
    where $\lVert \cdot \rvert\rvert_{A^p_\alpha}$ denotes the operator norm in $A^p_\alpha$.
\end{lm}

The operator $C_f^{\; n}\vcentcolon=C_{f^n}$ is bounded for all $n$ (both on $H^p$ and $A^p_\alpha$). In particular, by Lemmas \ref{lm:growth Hardy} and \ref{lm:growth Bergman}, the growth of the norms $\lVert C_f^{\,n}\rVert_{H^p}$ and $\lVert C_f^{\,n}\rVert_{A^p_\alpha}$ can be estimated by the quantity $1-\lvert f^n(0)\rvert$. But, due to \eqref{eq:euclidean rate} the quantities $d_{\D}(0,f^n(0))$ and $1-\lvert f^n(0)\rvert$ are equivalent. Thus, we can use Proposition \ref{prop:arosio bracci} of Arosio and Bracci to obtain the following:

\begin{cor}\label{cor:growth bounds}
    Let $f\vcentcolon \D\to\D$ be a non-elliptic map with Denjoy--Wolff point $\tau\in\partial\D$. Then
    \begin{enumerate}
        \item[\textup{(a)}] 
            $\lim\limits_{n\to+\infty}\dfrac{\log\lVert C_f^{\,n}\rVert_{H^p}}{n}=-\dfrac{\log f'(\tau)}{p}$, for all $p\ge1$;
        \item[\textup{(b)}] $\lim\limits_{n\to+\infty}\dfrac{\log\lVert C_f^{\,n}\rVert_{A^p_\alpha}}{n}=-\dfrac{(2+\alpha)\log f'(\tau)}{p}$, for all $p\ge1$ and all $\alpha>-1$.
    \end{enumerate}
\end{cor}

Note that if $f$ is parabolic, $f'(\tau)=1$ and hence both limits in Corollary \ref{cor:growth bounds} equal $0$. Replicating the arguments used to produce Corollary \ref{cor:growth bounds}, but using Theorem \ref{thm:main rates 2} instead of Proposition \ref{prop:arosio bracci}, we can obtain the precise estimate demonstrated in Corollary \ref{cor:operators}.

\subsection*{Acknowledgements}
We are grateful to the referees for their insightful comments and suggestions. For the first named author, the research project is implemented in the framework of H.F.R.I call ``3rd Call for H.F.R.I.'s Research Projects to Support Faculty Members \& Researchers” (H.F.R.I. Project Number: 24979).

\end{document}